\documentclass[12pt]{article}
\usepackage{amssymb}
\usepackage{graphicx}
\usepackage{latexsym}
\font\dynkfont=cmsy10 scaled\magstep4    \skewchar\dynkfont='60
\def\dynk{\textfont2=\dynkfont}
\def\hr#1,#2;{\dimen0=.4pt\advance\dimen0by-#2pt
              \vrule width#1pt height#2pt depth\dimen0}
\def\vr#1,#2;{\vrule height#1pt depth#2pt}
\def\blb#1#2#3#4#5
            {\hbox{\ifnum#2=0\hskip11.5pt
                   \else\ifnum#2=1\hr13,5;\hskip-1.5pt
                   \else\ifnum#2=2\hr13.5,6.5;\hskip-13.5pt
                                  \hr13.5,3.5;\hskip-2pt
                   \else\ifnum#2=3\hr13.7,8;\hskip-13.7pt
                                  \hr13,5;\hskip-13pt
                                  \hr13.7,2;\hskip-2.2pt\fi\fi\fi\fi
                   $#1$
                   \ifnum#4=0
                   \else\ifnum#4=1\hskip-9.2pt\vr22,-9;\hskip8.8pt
                   \else\ifnum#4=2\hskip-10.9pt\vr22,-8.75;\hskip3pt
                                  \vr22,-8.75;\hskip7.1pt
                   \else\ifnum#4=3\hskip-12.6pt\vr22,-8.5;\hskip3pt
                                  \vr22,-9;\hskip3pt
                                  \vr22,-8.5;\hskip5.4pt\fi\fi\fi\fi
                   \ifnum#5=0
                   \else\ifnum#5=1\hskip-9.2pt\vr1,12;\hskip8.8pt
                   \else\ifnum#5=2\hskip-10.9pt\vr1.25,12;\hskip3pt
                                  \vr1.25,12;\hskip7.1pt
                   \else\ifnum#5=3\hskip-12.6pt\vr1.5,12;\hskip3pt
                                  \vr1,12;\hskip3pt
                                  \vr1.5,12;\hskip5.4pt\fi\fi\fi\fi
                   \ifnum#3=0\hskip8pt
                   \else\ifnum#3=1\hskip-5pt\hr13,5;
                   \else\ifnum#3=2\hskip-5.5pt\hr13.5,6.5;
                                  \hskip-13.5pt\hr13.5,3.5;
                   \else\ifnum#3=3\hskip-5.7pt\hr13.7,8;
                                  \hskip-13pt\hr13,5;
                                  \hskip-13.7pt\hr13.7,2;\fi\fi\fi\fi
 
                   }}
\def\blob#1#2#3#4#5#6#7{\hbox
{$\displaystyle\mathop{\blb#1#2#3#4#5 }_{#6}\sp{#7}$}}
\def\up#1#2{\dimen1=33pt\multiply\dimen1by#1\hbox{\raise\dimen1\rlap{#2}}}
\def\uph#1#2{\dimen1=17.5pt\multiply\dimen1by#1\hbox{\raise\dimen1\rlap{#2}}}
\def\dn#1#2{\dimen1=33pt\multiply\dimen1by#1\hbox{\lower\dimen1\rlap{#2}}}
\def\dnh#1#2{\dimen1=17.5pt\multiply\dimen1by#1\hbox{\lower\dimen1\rlap{#2}}}
 
\def\rlbl#1{\kern-8pt\raise3pt\hbox{$\scriptstyle #1$}}
\def\llbl#1{\raise3pt\llap{\hbox{$\scriptstyle #1$\kern-8pt}}}
\def\elbl#1{\kern3pt\lower4.5pt\hbox{$\scriptstyle #1$}}
\def\lelbl#1{\rlap{\hbox{\kern-9pt\raise2.5pt\hbox{{$\scriptstyle #1$}}}}}

\def\whtd#1#2#3#4#5{\blob\circ#1#2#3#4{#5}{}}
\def\blkd#1#2#3#4#5{\blob\bullet#1#2#3#4{#5}{}}
\def\whtu#1#2#3#4#5{\blob\circ#1#2#3#4{}{#5}}
\def\blku#1#2#3#4#5{\blob\bullet#1#2#3#4{}{#5}}
\def\whtr#1#2#3#4#5{\blob\circ#1#2#3#4{}{}\rlbl{#5}}
\def\blkr#1#2#3#4#5{\blob\bullet#1#2#3#4{}{}\rlbl{#5}}

\def\rwng{\hbox{$\vbox{\offinterlineskip{
  \hbox{\phantom{}\kern6pt{$\circ$}}\kern-2.5pt\hbox{$\Biggr/$}\kern-0.5pt
  \hbox{\phantom{}\kern-5pt$\circ$}\kern-3.0pt\hbox{$\Biggr\backslash$}
  \kern-1.5pt\hbox{\phantom{}\kern6pt{$\circ$}} }}$}}
 
\def\lwng{\hbox{$\vbox{\offinterlineskip{ \hbox{$\circ$}
  \kern-3.0pt\hbox{\phantom{}\kern6.0pt{$\Biggr\backslash$}}
  \kern-0.5pt\hbox{\phantom{}\kern11pt{$\circ$}}\kern-3.5pt
  \hbox{\phantom{}\kern5.0pt {$\Biggr/$}}\kern-1.0pt\hbox{$\circ$} }}$}}

\def\drwng#1#2#3{\hbox{$\vcenter{ \offinterlineskip{
  \hbox{\phantom{}\kern7pt{$\circ^{\elbl{#3}}$}}
  \kern-2.5pt\hbox{$\Biggr/$}\kern-0.5pt
  \hbox{\phantom{}\kern-5pt$\circ^{ \elbl{#1}}$}
  \kern-3.0pt\hbox{$\Biggr\backslash$}
  \kern-1.5pt\hbox{\phantom{}\kern7pt{$\circ^{\elbl{#2}}$}}  } }$}}
\def\drwngt#1#2#3{\hbox{$\vcenter{ \offinterlineskip{
  \hbox{\phantom{}\kern7pt{$\bullet^{\elbl{#3}}$}}
  \kern-2.5pt\hbox{$\Biggr/$}\kern-0.5pt
  \hbox{\phantom{}\kern-5pt$\circ^{ \elbl{#1}}$}
  \kern-3.0pt\hbox{$\Biggr\backslash$}
  \kern-1.5pt\hbox{\phantom{}\kern7pt{$\bullet^{\elbl{#2}}$}}  } }$}}
 
\def\dlwng#1#2#3{\hbox{$\vcenter{\offinterlineskip{ \hbox{$\lelbl{#1}\circ$}
  \kern-3.0pt\hbox{\phantom{}\kern6.0pt{$\Biggr\backslash$}}
  \kern-0.5pt\hbox{\phantom{}\kern11pt{$\lelbl{#2}\circ$}}\kern-3.5pt
  \hbox{\phantom{}\kern5.0pt {$\Biggr/$}}\kern-1.0pt\hbox{$\lelbl{#3}\circ$}}}$}
 }

\def\rde#1#2#3{\hbox{\phantom{}\kern-4pt\hbox{$\vcenter{\offinterlineskip  \hbox
{
               \raise 4.5pt\hbox{\vrule height0.4pt width13pt depth0pt}
                \kern-1pt\vbox{ \hbox{\drwng{#1}{#2}{#3}}} }}$  }}  }
\def\rdet#1#2#3{\hbox{\phantom{}\kern-4pt\hbox{$\vcenter{\offinterlineskip \hbox
{
               \raise 4.5pt\hbox{\vrule height0.4pt width13pt depth0pt}
                \kern-1pt\vbox{ \hbox{\drwngt{#1}{#2}{#3}}} }}$  }}  }
 
\def\lde#1#2#3{\hbox{$\vcenter{\offinterlineskip  \hbox{
               \dlwng{#1}{#2}{#3}\kern-4.2pt\lower0.4pt\hbox{$\vcenter{\hrule 
                               width13pt}$}
               \kern-8pt\phantom{}   }}  $}}

\def\rwngb{\hbox{$\vbox{\offinterlineskip{
  \hbox{\phantom{}\kern6pt{$\bullet$}}\kern-2.5pt\hbox{$\Biggr/$}\kern-0.5pt
  \hbox{\phantom{}\kern-5pt$\bullet$}\kern-3.0pt\hbox{$\Biggr\backslash$}
  \kern-1.5pt\hbox{\phantom{}\kern6pt{$\bullet$}} }}$}}
 
\def\lwngb{\hbox{$\vbox{\offinterlineskip{ \hbox{$\bullet$}
  \kern-3.0pt\hbox{\phantom{}\kern6.0pt{$\Biggr\backslash$}}
  \kern-0.5pt\hbox{\phantom{}\kern11pt{$\bullet$}}\kern-3.5pt
  \hbox{\phantom{}\kern5.0pt {$\Biggr/$}}\kern-1.0pt\hbox{$\bullet$} }}$}}

\def\dbrwng#1#2#3{\hbox{$\vcenter{ \offinterlineskip{
  \hbox{\phantom{}\kern6pt{$\bullet^{\elbl{#3}}$}}
  \kern-2.5pt\hbox{$\Biggr/$}\kern-0.5pt
  \hbox{\phantom{}\kern-5pt$\bullet^{ \elbl{#1}}$}
  \kern-3.0pt\hbox{$\Biggr\backslash$}
  \kern-1.5pt\hbox{\phantom{}\kern6pt{$\bullet^{\elbl{#2}}$}}  } }$}}
 
\def\dblwng#1#2#3{\hbox{$\vcenter{\offinterlineskip{ \hbox{$\lelbl{#1}\bullet$}
  \kern-3.0pt\hbox{\phantom{}\kern6.0pt{$\Biggr\backslash$}}
  \kern-0.5pt\hbox{\phantom{}\kern11pt{$\lelbl{#2}\bullet$}}\kern-3.5pt
  \hbox{\phantom{}\kern5.0pt {$\Biggr/$}}\kern-1.0pt\hbox{$\lelbl{#3}\bullet$}}}
$} }

\def\rbde#1#2#3{\hbox{\phantom{}\kern-4pt\hbox{$\vcenter{\offinterlineskip  \hbo
x{
               \raise 4.5pt\hbox{\vrule height0.4pt width13pt depth0pt}
                \kern-1pt\vbox{ \hbox{\dbrwng{#1}{#2}{#3}}} }}$  }}  }
 
\def\lbde#1#2#3{\hbox{$\vcenter{\offinterlineskip  \hbox{
               \dblwng{#1}{#2}{#3}\kern-4.2pt\lower0.4pt\hbox{$\vcenter{\hrule w
idth13pt}$}
               \kern-8pt\phantom{}   }}  $}}

\def\ddgu#1.#2.{\dynk  \whtu0300{#1}\blku3000{#2}}
\def\ddgd#1.#2.{\dynk  \whtd0300{#1}\blkd3000{#2}}
 
\def\eddgiu#1.#2.#3.{\dynk \whtu0100{#1}\whtu1300{#2}\blku3000{#3}}
\def\eddgid#1.#2.#3.{\dynk \whtd0100{#1}\whtd1300{#2}\blkd3000{#3}}
\def\eddgiiu#1.#2.#3.{\dynk  \whtu0300{#1}\blku3100{#2}\blku1000{#3}}
\def\eddgiid#1.#2.#3.{\dynk  \whtd0300{#1}\blkd3100{#2}\blkd1000{#3}}

\def\ddfu#1.#2.#3.#4.{\dynk \whtu0100{#1}\whtu1200{#2}\blku2100{#3}\blku1000{#4}
}
\def\ddfd#1.#2.#3.#4.{\dynk \whtd0100{#1}\whtd1200{#2}\blkd2100{#3}\blkd1000{#4}
}
 
\def\eddfiu#1.#2.#3.#4.#5.{\dynk \whtu0100{#1}\whtu1100{#2}\whtu1200{#3}\blku210
0{#4}\blku1000{#5}}
\def\eddfid#1.#2.#3.#4.#5.{\dynk \whtd0100{#1}\whtd1100{#2}\whtd1200{#3}\blkd210
0{#4}\blkd1000{#5}}
\def\eddfiiu#1.#2.#3.#4.#5.{\dynk \whtu0100{#1}\whtu1200{#2}\blku2100{#3}\blku11
00{#4}\blku1000{#5}}
\def\eddfiid#1.#2.#3.#4.#5.{\dynk \whtd0100{#1}\whtd1200{#2}\blkd2100{#3}\blkd11
00{#4}\blkd1000{#5}}

\def\ddanu#1.#2.#3.#4.#5.{\dynk \whtu0100{#1}\whtu1100{#2}\whtu1100{#3}\cdots
                           \whtu1100{#4}\whtu1000{#5}}
\def\ddand#1.#2.#3.#4.#5.{\dynk \whtd0100{#1}\whtd1100{#2}\whtd1100{#3}\cdots
                           \whtd1100{#4}\whtd1000{#5}}

  \def\ddanm#1.#2.#3.#4.#5.{\dynk \whtu0100{#1}\whtu1100{#2}\whtu1100{#3}
                           \whtu1100{#4}\whtu1000{#5}} 
                           
    \def\ddnigel#1.#2.#3.{\dynk \whtu0100{#1}\whtu1100{#2}\whtu1100{#3}}
                                                       
       \def\ddanmm#1.#2.#3.#4.{\dynk \whtu0100{#1}\whtu1100{#2}
                           \whtu1100{#3}\whtu1000{#4}} 
                           
       \def\ddanmmm#1.#2.#3.{\dynk \whtu0100{#1}\whtu1100{#2}
                           \whtu1000{#3}}

 \def\ddans#1.#2.#3.#4.#5.{\dynk \blkd0100{#1}\whtd1100{#2}\blkd1100{#3}
                           \blkd1100{#4}\whtd1000{#5}}
\def\ddandto#1.#2.#3.#4.#5.{\dynk \blkd0100{#1}\whtd1100{#2}\blkd1100{#3}\cdots
                           \whtd1100{#4}\blkd1000{#5}}
 
\def\eddanu#1.#2.#3.#4.#5.{\dynk \whtu0100{#1}\whtu1100{#2}%
                           \up1{\whtr0000{#3}}\cdots\whtu1100{#4}\whtu1000{#5}}
\def\eddand#1.#2.#3.#4.#5.{\dynk \whtd0100{#1}\whtd1100{#2}%
                           \up1{\whtr0000{#3}}\cdots\whtd1100{#4}\whtd1000{#5}}
\def\eddaid#1.#2.{\dynk\whtd0400{#1}\hskip30pt\whtd4000{#2}}

\def\eddanid#1.#2.#3.#4.#5.{\dynk \whtd0200{#1}\whtd2100{#2}%
                           \whtd1100{#3}\cdots\whtd1200{#4}\blkd2000{#5}}
\def\eddaniu#1.#2.#3.#4.#5.{\dynk \whtu0200{#1}\whtu2100{#2}%
                           \whtu1100{#3}\cdots\whtu1200{#4}\blku2000{#5}}

\def\eddaniid#1.#2.#3.#4.#5.#6.{\hbox{$\vcenter{\hbox
         {\dynk\hbox{$ \lbde{#1}{#2}{#3}\blkd1100{#4}\cdots%
          \blkd1200{#5}\whtd2000{#6} $}} }$}}
\def\eddaniiu#1.#2.#3.#4.#5.#6.{\hbox{$\vcenter{\hbox
         {\dynk\hbox{$ \lbde{#1}{#2}{#3}\blku1100{#4}\cdots%
          \blku1200{#5}\whtu2000{#6} $}} }$}}

\def\eddaiiid#1.#2.{\dynk\blkd0400{#1}\hskip30pt\whtd4000{#2}}

\def\ddbnu#1.#2.#3.#4.#5.{\dynk \whtu0100{#1}\whtu1100{#2}\whtu1100{#3}\cdots
                           \whtu1200{#4}\blku2000{#5}}
\def\ddbnd#1.#2.#3.#4.#5.{\dynk \whtd0100{#1}\whtd1100{#2}\whtd1100{#3}\cdots
                           \whtd1200{#4}\blkd2000{#5}}
 
\def\eddbnu#1.#2.#3.#4.#5.#6.{\dynk \lde{#1}{#2}{#3}\whtu1100{#4}\cdots
                           \whtu1200{#5}\blku2000{#6}}
\def\eddbnd#1.#2.#3.#4.#5.#6.{\dynk \lde{#1}{#2}{#3}\whtd1100{#4}\cdots
                           \whtd1200{#5}\blkd2000{#6}}

\def\ddcnu#1.#2.#3.#4.#5.{\dynk \blku0100{#1}\blku1100{#2}\blku1100{#3}\cdots
                           \blku1200{#4}\whtu2000{#5}}
\def\ddcnd#1.#2.#3.#4.#5.{\dynk \blkd0100{#1}\blkd1100{#2}\blkd1100{#3}\cdots
                           \blkd1200{#4}\whtd2000{#5}}
 
\def\eddcnu#1.#2.#3.#4.#5.#6.{\dynk \whtu0200{#1}\blku2100{#2}\blku1100{#3}
       \blku1100{#4}\cdots
                           \blku1200{#5}\whtu2000{#6}}
\def\eddcnd#1.#2.#3.#4.#5.{\dynk \whtd0200{#1}\blkd2100{#2}\blkd1100{#3}
       \cdots \blkd1200{#4}\whtd2000{#5}}

\def\dddnu#1.#2.#3.#4.#5.#6.{\hbox{$\vcenter{\hbox
         {\dynk\hbox{$ \whtu0100{#1}\whtu1100{#2}\cdots%
          \whtu1100{#3}\rde{#4}{#5}{#6} $}}  }$}}
\def\dddnd#1.#2.#3.#4.#5.#6.{\hbox{$\vcenter{\hbox
         {\dynk\hbox{$ \whtd0100{#1}\whtd1100{#2}\cdots%
          \whtd1100{#3}\rde{#4}{#5}{#6} $}} }$}}
\def\dddndte#1.#2.#3.#4.#5.#6.{\hbox{$\vcenter{\hbox
         {\dynk\hbox{$ \blkd0100{#1}\whtd1100{#2}\cdots%
          \blkd1100{#3}\rdet{#4}{#5}{#6} $}} }$}}
\def\dddndto#1.#2.#3.#4.#5.#6.{\hbox{$\vcenter{\hbox
         {\dynk\hbox{$ \whtd0100{#1}\blkd1100{#2}\cdots%
          \blkd1100{#3}\rdet{#4}{#5}{#6} $}} }$}}
\def\dddiv#1.#2.#3.#4.{\hbox{$\vcenter{\hbox
         {\dynk\hbox{$ \whtu0100{#1}\rde{#2}{#3}{#4}
              $}}  }$}}

\def\edddnu#1.#2.#3.#4.#5.#6.#7.#8.{\hbox{$\vcenter{\hbox
         {\dynk\hbox{$ \lde{#1}{#2}{#3}\whtu1100{#4}\cdots%
          \whtu1100{#5}\rde{#6}{#7}{#8} $}}  }$}}
\def\edddnd#1.#2.#3.#4.#5.#6.#7.#8.{\hbox{$\vcenter{\hbox
         {\dynk\hbox{$ \lde{#1}{#2}{#3}\whtd1100{#4}\cdots%
          \whtd1100{#5}\rde{#6}{#7}{#8} $}} }$}}

\def\edddniid#1.#2.#3.#4.#5.{\hbox{$\vcenter{\hbox
         {\dynk\hbox{$ \blkd0200{#1}\whtd2100{#2}\whtd1100{#3}\cdots%
          \whtd1200{#4}\blkd2000{#5} $}} }$}}
\def\edddniiu#1.#2.#3.#4.#5.{\hbox{$\vcenter{\hbox
         {\dynk\hbox{$ \blku0200{#1}\whtu2100{#2}\whtu1100{#3}\cdots%
          \whtu1200{#4}\blku2000{#5} $}} }$}}

\def\ddei#1.#2.#3.#4.#5.#6.{\hbox{$\vcenter{\hbox
       {\dynk \whtd0100{#1}\whtd1100{#3}%
       \up1{\whtr0001{#2}}\whtd1110{#4}\whtd1100{#5}\whtd1000{#6}} }$}}
\def\ddeit#1.#2.#3.#4.#5.#6.{\hbox{$\vcenter{\hbox
       {\dynk \whtd0100{#1}\blkd1100{#3}%
       \up1{\blkr0001{#2}}\whtd1110{#4}\blkd1100{#5}\whtd1000{#6}} }$}}
 
\def\eddei#1.#2.#3.#4.#5.#6.#7.{\hbox{$\vcenter{\hbox
       {\dynk \whtd0100{#1}\whtd1100{#3}%
       \up1{\whtr0011{#2}}\up2{\whtr0001{#7}}\whtd1110{#4}\whtd1100{#5}%
       \whtd1000{#6}} }$}}

\def\ddeii#1.#2.#3.#4.#5.#6.#7.{\hbox{$\vcenter{\hbox
       {\dynk \whtd0100{#1}\whtd1100{#3}%
       \up1{\whtr0001{#2}}\whtd1110{#4}\whtd1100{#5}\whtd1100{#6}%
       \whtd1000{#7}} }$}}
\def\ddeiit#1.#2.#3.#4.#5.#6.#7.{\hbox{$\vcenter{\hbox
       {\dynk \whtd0100{#1}\blkd1100{#3}%
       \up1{\blkr0001{#2}}\whtd1110{#4}\blkd1100{#5}\whtd1100{#6}%
       \blkd1000{#7}} }$}}
 
\def\eddeii#1.#2.#3.#4.#5.#6.#7.#8.{\hbox{$\vcenter{\hbox
       {\dynk \whtd0100{#8}\whtd1100{#1}\whtd1100{#3}%
       \up1{\whtr0001{#2}}\whtd1110{#4}\whtd1100{#5}\whtd1100{#6}%
       \whtd1000{#7}} }$}}

\def\ddeiii#1.#2.#3.#4.#5.#6.#7.#8.{\hbox{$\vcenter{\hbox
       {\dynk \whtd0100{#1}\whtd1100{#3}%
       \up1{\whtr0001{#2}}\whtd1110{#4}\whtd1100{#5}\whtd1100{#6}%
       \whtd1100{#7}\whtd1000{#8}} }$}}
\def\ddeiiit#1.#2.#3.#4.#5.#6.#7.#8.{\hbox{$\vcenter{\hbox
       {\dynk \whtd0100{#1}\blkd1100{#3}%
       \up1{\blkr0001{#2}}\whtd1110{#4}\blkd1100{#5}\whtd1100{#6}%
       \blkd1100{#7}\whtd1000{#8}} }$}}
 
\def\eddeiii#1.#2.#3.#4.#5.#6.#7.#8.#9.{\hbox{$\vcenter{\hbox
       {\dynk \whtd0100{#1}\whtd1100{#3}%
       \up1{\whtr0001{#2}}\whtd1110{#4}\whtd1100{#5}\whtd1100{#6}%
       \whtd1100{#7}\whtd1100{#8}\whtd1000{#9}} }$}}

\def\part#1{\frac{\partial\phantom{q}}{\partial#1}}

\newenvironment{rmk}{\begin{trivlist}\item[]{\bf Remark:} }
{\end{trivlist}}
\newenvironment{rmks}{\begin{trivlist}\item[]{\bf Remarks:} }
{\end{trivlist}}
\newenvironment{ex}{\begin{trivlist}\item[]{\bf Example:} }
{\end{trivlist}}

\newenvironment{prf}{\begin{trivlist}\item[]{\bf Proof:} }
{\hfill $\Box$ \end{trivlist}}

\newtheorem{thm}{Theorem}

\newtheorem{prp}[thm]{Proposition}
\newtheorem{lemma}[thm]{Lemma}

\def\End{\mathop{\rm End}\nolimits}

\def\ev{\mathop{\rm ev}\nolimits}
\def\Pic{\mathop{\rm Pic}\nolimits}
\def\Hom{\mathop{\rm Hom}\nolimits}

\def\tr{\mathop{\rm tr}\nolimits}

\def\Sym{\mathop{\rm Sym}\nolimits}

\newcommand{\R}{\mathbf{R}}
\newcommand{\C}{\mathbf{C}}

\newcommand{\Z}{\mathbf{Z}}

\newcommand{\PP}{{\mathbf {\rm P}}}

\begin{document}
\title{Deformations of holomorphic Poisson manifolds}
 \author{Nigel Hitchin\\[5pt]}
 \maketitle
\centerline{{\it Subject classification}: {Primary }  32G05, 53D17, 53D18}
\section{Introduction}
One of the interesting by-products of generalized geometry is an unobstructedness theorem due to Goto for deformations of complex structure on a holomorphic Poisson manifold $M$  (Theorem 3.2 in  \cite{Goto}). This states that the contraction of the Poisson tensor $\sigma\in H^0(M,\Lambda^2T)$ with a K\"ahler form $\omega$ defines a class in $H^1(M,T)$ which can be integrated to a one-parameter family of deformations. The starting point for this paper is to prove this more generally, when $\omega$ is an arbitrary closed $(1,1)$-form and $M$ satisfies the $\partial\bar\partial$-lemma, without using the formalism of generalized geometry.  The argument, it turns out, is close to that of Bogomolov \cite{Bog} who used it in the symplectic case. 

This family of deformations parametrized by a single variable $t$ has some rather special properties. Each  deformation has an associated Poisson structure $\sigma_t$ and (if we view the deformation as a variation in the complex structure of a fixed $C^{\infty}$ manifold)   the subset on which  the tensor has a fixed rank is unchanged. Moreover, if $H^2(M,{\mathcal O})=0$ the Kodaira-Spencer class at $t$ is given by contraction  with the same cohomology class $[\omega]\in H^2(M,\C)\cong H^1(M,T^*)$. Furthermore, when the Poisson structure is generically symplectic, and hence defined by a closed meromorphic $2$-form, its cohomology class varies linearly in the direction $[\omega]$. We give examples of this deformation for some  classes of Poisson surfaces -- the projective plane $\PP_2$, the blow-up of $n$ points on a smooth cubic curve in $\PP_2$, a special rational elliptic surface related to Painlev\'e's first equation and compactifications of ALE spaces.

We then use this result to analyze deformations of the Hilbert scheme $M=S^{[n]}$ of points on a surface $S$. By a result of Fantechi  \cite{Fan} all deformations of $S^{[n]}$ for a surface of general type are induced from deformations of $S$, but her methods also show that for other surfaces there is a homomorphism $\rho: H^1(S^{[n]},T)\rightarrow H^0(S,K^*)$.  But this is  the space of Poisson structures on $S$, and   in \cite{Bot}  it is  shown that a Poisson structure on $S$ induces one on $S^{[n]}$. We show that if we take $\omega$ in the cohomology class of the exceptional divisor of  $S^{[n]}$ and apply the deformation theorem, we get  deformations of  $S^{[n]}$ whose Kodaira-Spencer classes  form a right inverse to $\rho$.

The most interesting case is where $S=\PP_2$ which, being rigid, means that all the deformations of the Hilbert scheme are obtained from Poisson structures  on $\PP_2$. We first investigate $\PP_2^{[2]}$ and describe it as a resolution of a cubic fourfold -- an intersection of the cubic $\det S=0$ in the $8$-dimensional projective space of $3\times 3$ matrices with a $5$-dimensional space. We then identify the higher-dimensional  case from the paper  \cite{Nev} of Nevins and Stafford. There, the authors produce a deformation of the Hilbert scheme $\PP_2^{[n]}$ as a moduli space of rank one modules over a non-commutative algebra. For the generic situation this is a Sklyanin algebra which depends on a smooth elliptic curve and a translation. We show that a generic deformation of the Hilbert scheme is of this form: moreover the one-parameter family of Poisson deformations constructed above preserve the modulus of the elliptic curve but change the translation linearly. 
\vskip .25cm
The author wishes to thank Daniel Huybrechts for useful communications and EPSRC for support. This work has been carried out in association with the ITGP network of the European Science Foundation. 

\section{Poisson geometry}
\subsection{Basics}
A Poisson structure on a  complex manifold $M$ is a holomorphic section $\sigma$ of $\Lambda^2T$ which satisfies an integrability condition. If $\sigma$ is non-degenerate then it defines a  holomorphic 2-form $\varphi$ and the integrability condition is $d\varphi=0$. In general, the condition is $[\sigma,\sigma]=0\in H^0(M,\Lambda^3T)$ using the Schouten-Nijenhuis bracket, and in particular this condition is automatic in two dimensions.  A more useful  description of integrability is to consider $\sigma$ as a homomorphism $\sigma:T^*\rightarrow T$, and to  take a local holomorphic function $f$ and define the corresponding Hamiltonian vector field by $\sigma(df)$. Then integrability is equivalent to 
$$[\sigma(df),\sigma(dg)]=\sigma(d \{f,g\})$$
where $\{f,g\}$ is the Poisson bracket $\sigma(df)g=-\sigma(dg)f$. Integrability means that the Poisson bracket satisfies the Jacobi identity.
In local coordinates we shall write,  using the summation convention, 
$$\sigma=\sigma^{ij}\frac{\partial}{\partial z_i}\wedge\frac{\partial}{\partial z_j}.$$

 A good survey of Poisson geometry in the algebro-geometric context is \cite{P1}. We shall be entirely concerned with compact Poisson manifolds in this paper. Compact Poisson surfaces are classified in \cite{BM} but it is not so easy to find examples in higher dimensions. One general case is give by the canonical symplectic structure on the cotangent bundle $T^*$ of a complex manifold. This extends to a Poisson structure on its compactification $\PP(T^*\oplus {\mathcal O})$. In \cite{P2} some concrete Poisson structures on projective spaces and Grassmannians are defined as moduli spaces of chains of bundles on an elliptic curve. It is a general result \cite{Bot0} that moduli spaces of sheaves on a Poisson surface are Poisson. The special case of Hilbert schemes \cite{Bot} yields smooth compact examples, and we shall return to this case in Section  \ref{HS}. 

\subsection{Deformations}\label{def}
For a deformation of a complex structure, one looks for  a global $\phi\in \Omega^{0,1}(T)$ which satisfies the equation
\begin{equation}
\bar\partial\phi +\frac{1}{2}[\phi,\phi]=0
\label{MC}
\end{equation}
where the bracket is the Lie bracket on  vector fields together with  exterior product on $(0,1)$-forms.  
Given such a $\phi=\phi_{\bar i}d\bar z_i$, then 
$$\left[\frac{\partial}{\partial \bar z_i}+\phi_{\bar i},\frac{\partial}{\partial \bar z_j}+\phi_{\bar j}\right]=0.$$
If the matrix $(1-\phi\bar\phi)_{ij}=\delta^i_j-\phi^i_{\bar k}\bar\phi^{\bar k}_j$ is invertible, then  it follows that we have a complex structure whose $(0,1)$ vector fields are spanned by  the commuting vector fields 
$$\bar X_i=\frac{\partial}{\partial \bar z_i}+\phi_{\bar i}=\frac{\partial}{\partial \bar z_i}+\phi^j_{\bar i}\frac{\partial}{\partial z_j}.$$
To obtain a deformation one tries to solve  (\ref{MC})  term-by-term for a series $\phi(t)=t\phi_1+t^2\phi_2+\dots$.  Having done this, for small enough $t$  elliptic estimates prove convergence and we also have invertibility of $(1-\phi\bar\phi)$.  For such a deformation  $\bar\partial\phi_1=0$ and $\phi_1$ represents the Kodaira-Spencer class of the deformation in $H^1(M,T)$. Conversely given such a class one may ask if a deformation exists in that direction. 

On a Poisson manifold there is a natural way to obtain such classes. Let $\omega$ be a  $(1,1)$-form with $\bar\partial\omega=0$. Then applying $\sigma:T^*\rightarrow T$ we obtain 
$$\phi_1=\sigma(\omega)\in \Omega^{0,1}(T)$$
which is  $\bar\partial$-closed since $\sigma$ is holomorphic. Cohomologically this is just the natural contraction map $H^0(M,\Lambda^2T)\otimes H^1(M,T^*)\rightarrow  H^1(M,T)$.

\begin{rmks}

\noindent 1. This process sometimes gives a trivial class in $H^1(M,T)$, for example if $[\omega]$ is a multiple of the first Chern class of $M$. In this case $-c_1$ is the Atiyah class of the canonical bundle $K$ -- the obstruction to the existence of a holomorphic connection. When we apply $\sigma$ we get the obstruction to the existence of a holomorphic first order differential operator   ${\mathcal D}:K\rightarrow  K\otimes T$ whose symbol is $\sigma$. But on $K$ there exists  such an operator characterized by 
${\mathcal D}_{df}s={\mathcal L}_{\sigma(df)}s.$ Indeed, if $X=\sigma(df)=a^i\partial/\partial z_i$ we have 
$${\mathcal L}_X(dz_1\wedge\dots\wedge dz_n)=\frac{\partial a^i}{\partial z_i} (dz_1\wedge\dots\wedge dz_n)$$
and, since $\sigma^{ij}$ is skew-symmetric,
$$\frac{\partial a^i}{\partial z_i}=\frac{\partial}{\partial z_i}\left(\sigma^{ij}\frac{\partial f}{\partial z_j}\right)=\frac{\partial \sigma^{ij}}{\partial z_i}\frac{\partial f}{\partial z_j}$$
which is linear in  the first  derivative of $f$. Hence we can define the derivative ${\mathcal D}_{\alpha}s$ for any $(1,0)$-form $\alpha$, not just $df$.

The operator ${\mathcal D}$ also satisfies the ``zero curvature" condition ${\mathcal D}^2s=0\in {\mathcal O}(K\otimes \Lambda^2 T)$. This makes $K$ a {\it Poisson module}. 

\noindent 2. Another case is the first Chern class of the line bundle  defined by an irreducible  {\it component} $C$ of the anticanonical divisor  of $\sigma$ on a surface. If $L$ is the corresponding line bundle and $s$ the section vanishing on $C$, then the Atiyah class  is $\delta(ds)\in H^1(M,T^*)$ where $\delta$ is the coboundary map in the long exact sequence of 
$$0\rightarrow T^*\stackrel{s}\rightarrow  LT^*\rightarrow LT^*\vert_C\rightarrow 0.$$
Since $\sigma$ vanishes on $C$, $\sigma(ds)=0$ and hence $\sigma\delta(ds)=\delta\sigma(ds)=0$.
\end{rmks}

We shall prove in the theorem below that all such Kodaira-Spencer classes can be integrated to a finite deformation if $M$ satisfies the $\bar\partial\partial$-lemma. 

\begin{ex} An example (though not our principal concern here) of such a deformation is the twistor deformation of a hyperk\"ahler manifold. We have  complex structures $I,J,K$ satisfying  the algebraic relations of quaternions and corresponding K\"ahler forms $\omega_1,\omega_2,\omega_3$. With respect to $I$, $\varphi=\omega_2+i\omega_3$ is a holomorphic symplectic structure and hence defines a holomorphic Poisson structure. The closed $(1,1)$ form $\omega_1$ defines a Kodaira-Spencer class $\sigma(\omega_1)$ which integrates to the family of complex structures $\cos t I+\sin t K$. 
\end{ex}

\begin{thm} \label{deform}Let $(M,\sigma)$ be a holomorphic Poisson manifold which satisfies the $\partial\bar\partial$-lemma. Then any class $\sigma([\omega])\in H^1(M,T)$ for $[\omega]\in H^1(M,T^*)$ is tangent to a deformation of complex structure.
\end{thm}
\begin{prf}
Since the $\partial\bar\partial$-lemma holds we can represent the class $[\omega]$ by a closed $(1,1)$-form $\omega$ and we then need to solve equation (\ref{MC}) term-by-term for $\phi(t)=t\phi_1+t^2\phi_2+\dots$ where $\phi_1=\sigma(\omega)$ or 
$$\phi_1=\sigma^{i j}\omega_{j \bar k}\frac{\partial}{\partial z_i}\,d\bar z_k.$$
The  coefficient of $t^2$ requires a solution for $\phi_2$ of 
\begin{equation}
\bar\partial \phi_2+\frac{1}{2}[\phi_1,\phi_1]=0. 
\label{ob}
\end{equation}
 Locally we  write  $\omega=\partial \bar\partial h$. Set  $\bar\partial h=\alpha=a_{\bar k}d\bar z_k$ and then 
$$\omega_{ j \bar k}=\frac{\partial a_{\bar k} }{ \partial  z_{j}}$$
so that 
$$\phi_1=\sigma^{i j}\omega_{j \bar k}\frac{\partial}{\partial z_i}\,d\bar z_k=\sigma^{i j}\frac{\partial a_{\bar k} }{ \partial  z_{j}}\frac{\partial}{\partial z_i}\,d\bar z_k=\sigma(\partial a_{\bar k})d\bar z_k.$$
Therefore 
$$[\phi_1,\phi_1]=\sigma(\partial\{a_{\bar j},a_{\bar k}\})d\bar z_j\wedge d\bar z_{k}$$
using the integrability property $[\sigma(\partial f),\sigma(\partial g)]=\sigma(\partial \{f,g\})$ of $\sigma$. 

The Poisson bracket expression $\{a_{\bar i},a_{\bar j}\}d\bar z_i\wedge d\bar z_{j}$ looks local  but $\{f,g\}=\sigma(\partial f,\partial g)$ so it is 
$\sigma^{k\ell}\omega_{k\bar i}\omega_{\ell\bar j}d\bar z_i\wedge d\bar z_{j}$
or  ${\sigma}(\omega^2)$. Since  $\sigma$ is holomorphic and $\omega$ is $\bar\partial$-closed,  ${\sigma}(\omega^2)$ is  also $\bar\partial$-closed.

Thus $\partial ({\sigma}(\omega^2))$ is $\bar\partial$-closed and $\partial$-exact and so, by the $\partial\bar\partial$-lemma 
$$\partial ({\sigma}(\omega^2))=\bar\partial\partial\beta$$
for some $(0,1)$-form $\beta$. It follows that 
$$[\phi_1,\phi_1]=\sigma(\partial ({\sigma}(\omega^2)))=\sigma(\bar\partial\partial\beta)=\bar \partial( \sigma(\partial\beta))$$
and we take $\phi_2=-\sigma(\partial\beta)/2$ to solve Equation (\ref{ob}). Note that $\phi_2$ has the same form $\sigma(\omega)$ as $\phi_1$ but now $\omega$ is replaced by $-\partial\beta/2$, which  is $\partial$-exact. Write $\beta_2=-\beta/2$ and $\beta_1=\alpha=\bar\partial h$ then $\phi_k=\sigma(\partial\beta_k)$ for $k=1,2$. (For convenience we shall keep the notation $\beta_1$ even though it is only locally defined. In most of what follows it appears in a Poisson bracket which factors through $\partial\beta_1$ which {\it is} globally defined). 

Inductively, suppose that 
$$\bar\partial\partial \beta_k=-\frac{1}{2}\partial(\{\beta_1,\beta_{k-1}\}+\{\beta_2,\beta_{k-2}\}+\dots+\{\beta_{k-1},\beta_1\}).$$
for $k<n$. Now consider
$$ \gamma_n=\{\beta_1,\beta_{n-1}\}+\{\beta_2,\beta_{n-2}\}+\dots+\{\beta_{n-1},\beta_1\}.$$
Each term $\{\beta_k,\beta_{n-k}\}$ can be written $\sigma(\partial\beta_k\partial\beta_{n-k})$ and so
$$\bar\partial\{\beta_k,\beta_{n-k}\}=\sigma(\bar\partial\partial\beta_k\partial\beta_{n-k})-\sigma(\partial\beta_k\bar\partial\partial\beta_{n-k})$$
and by the inductive assumption this is
$$-\frac{1}{2}\{\{\beta_1,\beta_{k-1}\}+\dots+\{\beta_{k-1},\beta_1\},\beta_{n-k}\}+ \frac{1}{2}\{\beta_{k},\{\beta_1,\beta_{n-k-1}\}+\dots+\{\beta_{n-k-1},\beta_1\}\}. $$
Summing over $k$ this is 
$$\sum_{i+j+k=n}\{\beta_i,\{\beta_j,\beta_k\}\} $$
which vanishes by the Jacobi identity for the Poisson bracket, so $\bar\partial\gamma_n=0$.

Hence $\partial\gamma_n=\partial(\{\beta_1,\beta_{n-1}\}+\{\beta_2,\beta_{n-2}\}+\dots+\{\beta_{n-1},\beta_1\})$ is $\bar\partial$-closed and $\partial$-exact, so by the $\partial\bar\partial$-lemma can be written as $\bar\partial\partial(-2\beta_n)$ for some $\beta_n$, completing the induction.

Now define $\phi_k=\sigma(\partial\beta_k)$, then
\begin{eqnarray*}
\bar\partial\phi_k&=&-\frac{1}{2}\sigma(\partial(\{\beta_1,\beta_{k-1}\}+\{\beta_2,\beta_{k-2}\}+\dots+\{\beta_{k-1},\beta_1\})\\
&=&-\frac{1}{2}([\sigma\partial\beta_1,\sigma\partial\beta_{k-1}]+[\sigma\partial\beta_2,\sigma\partial\beta_{k-2}]+\dots+[\sigma\partial\beta_{k-1},\sigma\partial\beta_1])\\
&=&-\frac{1}{2}([\phi_1,\phi_{k-1}]+[\phi_2,\phi_{k-2}]+\dots+[\phi_{k-1},\phi_1])\\
\end{eqnarray*}
as required for the deformation.
\end{prf}

\begin{rmk} Note that if $H^2(M,{\mathcal O})=0$, then the $\bar\partial$-closed $(0,2)$ form $\gamma_n$ in the proof is $\bar\partial$-exact and we can define $\beta_n$ in the induction by  $\gamma_n=2\bar\partial \beta_n$ without using the $\partial\bar\partial$-lemma. If further $H^1(M,{\mathcal O})=0$, then  $\beta_n$ is unique modulo $\bar\partial f_n$ which generates a time-dependent Hamiltonian vector field. Under these circumstances the deformation is uniquely determined up to Poisson diffeomorphism. By contrast, in the hyperk\"ahler case, where $H^2(M,{\mathcal O})\ne 0$, the Kodaira-Spencer class is tangential to many one-parameter families of deformations, the hyperk\"ahler rotation being just one of them.
\end{rmk}

The theorem shows that globally we have a deformation given by 
$$\phi=\sigma(t\omega+\partial(t^2\beta_2+t^3\beta_3+\dots))$$
or locally by $\phi=\sigma(\partial\beta)$ for $\beta=t\beta_1+t^2\beta_2+\dots$ where $\beta$ is a $(0,1)$-form with respect to the initial complex structure. From the inductive  part of the proof we have 
$$\bar\partial\partial\beta+\frac{1}{2}\partial\{\beta,\beta\}=0.$$

 To describe the complex structure at $t=a$  more concretely, we take as above a local basis $X_1,\dots,X_n$ of $(1,0)$ vector fields 
$$X_i=\frac{\partial}{\partial z_i}+\bar\phi^{\bar j}_{ i}\frac{\partial}{\partial \bar z_j}$$
and the corresponding basis of $(1,0)$-forms $\xi_1,\dots,\xi_n$ 
$$\xi_i=(1-\bar\phi\phi)^{-1}_{ji}(dz_j-\phi^j_{\bar k}d\bar z_k).$$
Then for any function $f$ 
$$\bar\partial_a f=(\bar X_if)\bar\xi_i=\left(\frac{\partial f}{\partial \bar z_i}+\phi^{ j}_{ \bar i}\frac{\partial f}{\partial  z_j}\right)\bar\xi_i.$$
(In particular note that $\bar\partial_a\bar z_i=\bar\xi_i$ and so $\bar\partial_a\bar\xi_i=0$.) Using the Poisson  bracket we may write this also as 
\begin{equation}
\bar\partial_a f=\left(\frac{\partial f}{\partial \bar z_i}+\{\beta_{\bar i},f\}\right)\bar\xi_i.
\label{dbar}
\end{equation}

\subsection{The deformed Poisson structure}
The above theorem gives us a one-parameter  family of deformations of $M$ as a complex manifold. We shall see firstly that each such deformation is also a holomorphic Poisson manifold. If $f,g$ are local holomorphic functions with respect to the complex structure at $t=a$, we define
$$\sigma_a(df,dg)=\sigma(\partial f,\partial g).$$
In fact because $\sigma$ is a bivector of type $(2,0)$ its interior product with a $(0,1)$-form vanishes so we could as well write $\sigma_a(df,dg)=\sigma(df,dg).$ Using our local basis of $(1,0)$-forms we have 
\begin{equation}
\sigma_a(\xi_i,\xi_j)=\sigma(\xi_i,\xi_j)=(1-\bar\phi\phi)^{-1}_{ki}\sigma(dz_k,dz_{\ell})(1-\bar\phi\phi)^{-1}_{\ell j}
\label{sigmaa}
\end{equation}
\begin{prp} \label{holp} The bivector field $\sigma_a$
 is a holomorphic Poisson structure.
\end{prp}
\begin{prf}
Let $f,g$ be holomorphic with respect to the complex structure at $t=a$. Then $\bar X_k f=0=\bar X_k g$ and so 
$${\mathcal L}_{\bar X_k}(df)=0={\mathcal L}_{\bar X_k}(dg).$$
Now  $\bar X_k={\partial}/{\partial \bar z_k}+\phi_{\bar k}$ and   $\phi_{\bar k}=\sigma(\partial \beta_{\bar k})$, which  is a  complex Hamiltonian vector field and so the Lie derivative of $\sigma$ vanishes. (Notice that the only derivatives involved in the proof of this statement are with respect to $z_i$ and not $\bar z_i$ so it is immaterial whether  $\beta_{\bar k}$ is holomorphic or not.)
But $\sigma$ is also holomorphic in the original complex  structure so its Lie derivative by $\partial/\partial \bar z_k$ is  also zero. It follows that $${\mathcal L}_{\bar X_k}\sigma=0$$ 
and hence ${\bar X_k}\sigma(df,dg)=0$ for all holomorphic $f,g$ and all $k$, i.e. $\sigma_a$ is holomorphic.

The integrability condition for the Poisson structure is the Jacobi identity
$$\{f,\{g,h\}\}+\{h,\{f,g\}\}+\{g,\{h,f\}\}=0.$$
 
Now $\{f,\{g,h\}\}=\sigma_a(df,d(\sigma_a(dg,dh)))$ and for  local holomorphic functions $f,g,h$ $\sigma_a(dg,dh)=\sigma(dg,dh)$ which we have just shown is  also holomorphic and thus
$$\{f,\{g,h\}\} =\sigma(df,d(\sigma(dg,dh)))=\sigma(\partial f,\partial (\sigma(\partial g,\partial h)))$$
so integrability follows from the integrability of $\sigma$.
\end{prf}

\begin{rmks} 

\noindent 1. Note from (\ref{sigmaa}) that the subset of $M$ on which the rank of $\sigma=2k$ is unchanged under deformation. In particular this applies to the set where $\sigma=0$.  But here $\phi=\sigma(\partial\beta)$  itself vanishes and so not only is the zero set of the Poisson structure unchanged, but its holomorphic structure too. 

\noindent 2. Producing a new Poisson structure by restricting the old one to a new set of $(1,0)$-forms  is something which also occurs in a hyperk\"ahler manifold. The real  K\"ahler forms $\omega_i$ together with the metric define real Poisson structures $\sigma_i$ and there is a natural holomorphic Poisson structure $\sigma_{\zeta}=(\sigma_2+i\sigma_3)+2i\zeta\sigma_1+\zeta^2(\sigma_2-i\sigma_3)$ for the complex structure parametrized by $\zeta\in \C$ in the twistor family. To see what the restriction of this is to the complex structure at $\zeta$  it suffices to consider the flat case of $\C^{2n}$ with complex coordinates $z_1,\dots,z_n,w_1,\dots,w_n$.  Here  
$$\sigma_{\zeta}=\left(\frac{\partial}{\partial z_i}+\zeta \frac{\partial}{\partial \bar w_i}\right)\wedge\left(\frac{\partial}{\partial w_i}-\zeta \frac{\partial}{\partial \bar z_i}\right)=X_i\wedge Y_i.$$
The $(1,0)$ forms for the complex structure at $\zeta$ are spanned by the dual basis $\xi_i,\eta_j$ to $X_i,Y_j$ given by 
$$\xi_i=\frac{1}{1+\vert \zeta\vert^2}(dz_i+\bar\zeta d\bar w_i)\qquad\eta_j=
\frac{1}{1+\vert \zeta\vert^2}(dw_j-\bar\zeta d\bar z_j).$$
Then
$$\sigma_0(\xi_i,\eta_j)=\frac{1}{(1+\vert \zeta\vert^2)^2}\delta_{ij}=\frac{1}{(1+\vert \zeta\vert^2)^2}\sigma_{\zeta}(\xi_i,\eta_j)$$
and so $\sigma_0$ restricted to the new $(1,0)$ forms is a multiple of the hyperk\"ahler Poisson structure. 
\end{rmks}

\subsection{The Kodaira-Spencer class}
Since the deformation at $t=a$ has a natural Poisson structure $\sigma_a$, we can ask whether the Kodaira-Spencer class at $a$ is again defined by contraction of $\sigma_a$ with a closed $(1,1)$-form. To do this we   work  in the complex structure at $a$ using the  local basis $\xi_i$ of $(1,0)$-forms. Note from (\ref{dbar}) that $\bar\partial_a \bar z_i=\bar\xi_i$. 

Consider the local $(0,1)$-form $\beta(t,z,\bar z)=\beta_{\bar i}d\bar z_i$ used in the construction of the deformation. It is defined on an open set $U$ and we denote its $t$-derivative at $t=a$ by $\dot\beta$.    We write
$$\gamma_U= \dot\beta_{\bar i}\bar\xi_i$$
which is a  $(0,1)$-form in the complex structure at $a$.  We first prove the following:

\begin{prp} Suppose $H^2(M,{\mathcal O})=0$, then there is a well-defined closed $(1,1)$-form $\omega_a$ on $M$ such that $\partial_a\gamma_U=\omega_a\vert_U$.
\end{prp}

\begin{rmk} By semi-continuity, if $H^2(M,{\mathcal O})$ vanishes at $t=0$ then it also does for all small enough $t$. The same holds for the $\partial\bar\partial$-lemma and so all degree two cohomology classes are represented by closed $(1,1)$-forms for these  deformations.
\end{rmk}

\begin{prf}
From the formula (\ref{dbar}) for $\bar\partial_af$ and $\bar\partial_a\bar\xi_i=0$ we have 
\begin{equation}
\bar\partial_a \gamma_U=\left(\frac{\partial \dot\beta_{\bar i}}{\partial \bar z_j}+\{\beta_{\bar j},\dot\beta_{\bar i}\}\right)\bar\xi_j\bar\xi_i.
\label{gamma}
\end{equation}
If $H^2(M,{\mathcal O})=0$ then as remarked above in the induction in the proof we can take 
$$ \bar\partial\beta_n=\frac{1}{2}(\{\beta_1,\beta_{n-1}\}+\{\beta_2,\beta_{n-2}\}+\dots+\{\beta_{n-1},\beta_1\}).$$
Together with $\bar\partial\beta_1=0$ this gives 
$$\bar\partial\beta=\frac{1}{2}\{\beta,\beta\}$$
and differentiating with respect to $t$, $\bar\partial\dot\beta=\{\dot\beta,\beta\}$. Putting this in (\ref{gamma}) gives $\bar\partial_a\gamma_U=0$. Thus $\partial_a\gamma_U$ is a locally defined closed $(1,1)$-form. 

Let $\omega=\partial\bar\partial h_U$ on $U$. Now because all terms but the first $\beta_1=\bar\partial h_U$  in the expansion of $\beta$ are globally defined, on  $U\cap V$ we have 
\begin{equation}
\gamma_V-\gamma_U=\frac{\partial (h_V-h_U)}{\partial \bar z_i}\bar\xi_i.
\label{bound1}
\end{equation}
 Since $\omega$ is globally defined, on $U\cap V$ we have $\bar\partial(\partial h_V-\partial h_U)=0$ so $(\partial h_V-\partial h_U)$ is a $1$-cocycle with values in the sheaf $d{\mathcal O}$. From the  exact sequence of sheaves
$$0\rightarrow d{\mathcal O}\rightarrow \Omega^1\rightarrow d\Omega^1\rightarrow 0$$
it defines the class $[\omega]\in H^1(M,T^*)$. Now consider the exact sequence
$$0\rightarrow \C\rightarrow {\mathcal O}\rightarrow  d{\mathcal O}\rightarrow 0$$
and write $(\partial h_V-\partial h_U)=\partial g_{UV}$ where $g_{UV}$ is holomorphic on $U\cap V$. We have $\partial (h_V- h_U- g_{UV})=0$ and so $h_V- h_U- g_{UV}=g'_{UV}$ is an antiholomorphic function. But then $\bar\partial_a g'_{UV}$ can be written in terms of  $\bar\partial g'_{UV}$  and we obtain 
$$\bar\partial_a g'_{UV}=\frac{\partial g'_{UV}}{\partial \bar z_i}\bar\xi_i=\frac{\partial}{\partial \bar z_i}(h_V- h_U- g_{UV})\bar\xi_i=\frac{\partial}{\partial \bar z_i}(h_V- h_U)\bar\xi_i.$$
So from (\ref{bound1})  
$\gamma_V-\gamma_U=\bar\partial_ag'_{UV}.$ 
It follows that 
$$\partial_a\gamma_V-\partial_a\gamma_U=\partial_a\bar\partial_ag'_{UV}=d\bar\partial_ag'_{UV}$$
and since $g'_{UV}$ is antiholomorphic $\bar\partial_ag'_{UV}=\bar\partial g'_{UV}$ so
$$d\bar\partial_ag'_{UV}=d\bar\partial g'_{UV}=-\bar\partial\partial g'_{UV}=0.$$ 

Thus $\omega_a=\partial_a\gamma_U$ is a globally defined closed $(1,1)$-form in the complex structure at $t=a$. 
\end{prf}

\begin{prp} \label{coh}The closed forms $\omega$ and $\omega_a$ represent the same class in $H^2(M,\C)$.
\end{prp}

\begin{prf} This is a continuation of the \v{C}ech argument. With $\omega=\partial\bar\partial h_U$ we obtained  holomorphic functions $g_{UV}$ such that $(\partial h_V-\partial h_U)=\partial g_{UV}$. Then on $U\cap V\cap W$ we have $\partial (g_{UV}+g_{VW}+g_{WU})=0$ and so a constant $c_{UVW}=g_{UV}+g_{VW}+g_{WU}$ which is a $\C$-valued 2-cocycle representing the cohomology class of $\omega$ in $H^2(M,\C)$. Interchanging the roles of $\partial$ and $\bar\partial$ we have $\omega=-\bar\partial\partial h_U$ and we obtain the class of $-c_{UVW}$.

Now if we locally write $\gamma_U=\bar\partial_ak_U$ then $\omega_a=\partial_a\bar\partial_ak_U$ and 
$$\gamma_V-\gamma_U=\bar\partial_a(k_V-k_U)=\bar\partial_ag'_{UV}$$ 
so $c'_{UVW}=g'_{UV}+g'_{VW}+g'_{WU}$ defines the negative of the class of $\omega_a$. But we saw above that
$h_V- h_U- g_{UV}=g'_{UV}$ and so $c_{UVW}=-c'_{UVW}$ hence the two forms have the same cohomology class.
\end{prf}

\begin{prp} \label{KS} If $H^2(M,{\mathcal O})=0$, then the Kodaira-Spencer class of the deformation at $t=a$ is defined by $\sigma_a(\omega_a)$.
\end{prp}

\begin{prf}
The Kodaira-Spencer class at $t=a$ is obtained by taking the $(1,0)$ part of the $t$-derivative of the $(0,1)$ vector fields. Using the standard local basis, this is given by
$$\xi_j\left(\frac{d\bar X_i}{dt}\right)X_j\bar \xi_i=(1-\bar\phi\phi)^{-1}_{kj}\dot\phi^k_{\bar i}X_j\bar \xi_i=(1-\bar\phi\phi)^{-1}_{kj}\sigma^{k \ell}\frac{\partial\dot\beta_{\bar i}}{\partial z_{\ell}} X_j\bar \xi_i.$$
Now in this basis
 \begin{equation}
\omega_a=\partial_a( \dot\beta_{\bar i}\bar\xi_i)=(\partial_a \dot\beta_{\bar i})\bar\xi_i+ \dot\beta_{\bar i}\partial_a\bar\xi_i.
\label{omegabasis}
\end{equation}
We now need a lemma:

\begin{lemma} The $(0,1)$-form $\sigma_a(\partial_a\bar\xi_i,\xi_j)$ vanishes.
\end{lemma}
\begin{prf} We saw in Proposition \ref{holp} that ${\mathcal L}_{\bar X_k}\sigma=0$ and ${\mathcal L}_{\bar X_k}(df)=0$ for a holomorphic function $f$ (with respect to the complex structure at $t=a$). Since $\sigma_a$ is defined as $\sigma$ restricted to the derivatives of holomorphic functions it follows that 
${\mathcal L}_{\bar X_k}\sigma_a=0.$

Now since $i_{\bar X_k}\bar\xi_i=\delta_{ik}$, we have ${\mathcal L}_{\bar X_k}\bar \xi_i=i_{\bar X_k}d\bar\xi_i+d(i_{\bar X_k}\bar\xi_i)=i_{\bar X_k}d\bar\xi_i$. Furthermore, since $\bar\partial_a\bar\xi_i=0$, $i_{\bar X_k}d\bar\xi_i=i_{\bar X_k}\partial_a\bar\xi_i$. Hence 
$$i_{\bar X_k} \sigma_a(\partial_a\bar\xi_i,\xi_j)=\sigma_a({\mathcal L}_{\bar X_k}\bar \xi_i,\xi_j).$$

But $\sigma_a$ is of type $(2,0)$ so  $\sigma_a(\bar\xi_i, \eta)=0$ for all $\eta$ hence 
$$\sigma_a({\mathcal L}_{\bar X_k}\bar \xi_i,\xi_j)=-({\mathcal L}_{\bar X_k}\sigma_a)(\bar \xi_i,\xi_j)+{\bar X_k}(\sigma_a(\bar \xi_i,\xi_j))-\sigma_a(\bar \xi_i,{\mathcal L}_{\bar X_k}\xi_j)=0.$$

Thus $i_{\bar X_k} \sigma_a(\partial_a\bar\xi_i,\xi_j)=0$ for all $k$, proving the lemma.
\end{prf}

Using the lemma we can write, using equation (\ref{omegabasis}), 
$$\sigma_a(\omega_a)=\sigma_a(\partial_a \dot\beta_{\bar i},\xi_j)X_j\bar\xi_i.$$
But $\sigma_a$ evaluated on $(1,0)$-forms is just the restriction of $\sigma$, which is of type $(2,0)$ in the original complex structure.  This annihilates the $(0,1)$ components of $\partial_a \dot\beta_{\bar i}$ and $\xi_j$ so  we get 
$$ \sigma(\partial_a \dot\beta_{\bar i},\xi_j)X_j\bar\xi_i = (1-\bar\phi\phi)^{-1}_{kj}\sigma^{k \ell}\frac{\partial\dot\beta_{\bar i}}{\partial z_{\ell}} X_j\bar \xi_i $$
thereby proving the Proposition.
\end{prf}

We see here that there is nothing special about $t=0$ in this deformation family if we use the Poisson structure $\sigma_a$: each Kodaira-Spencer class in the deformation is  given  by contraction of the Poisson tensor with a closed $(1,1)$ form in the same cohomology class.

\subsection{Periods}

In many examples (in particular in two dimensions) the Poisson tensor $\sigma$ is generically non-degenerate and so its inverse $\varphi$ defines a closed meromorphic $2$-form with a pole along an anticanonical divisor $D$. On $M\backslash D$ this form is regular and so has a cohomology class in $H^2(M\backslash D,\C)$. Since the Poisson deformation preserves the subset $D$ where the Poisson tensor drops rank, there is a corresponding $2$-form $\varphi_a$ on $M\backslash D$ and we can ask how the periods vary in the deformation. We have 

\begin{prp} \label{pi} Let $M$, with $H^2(M,{\mathcal O})=0$, be a holomorphic Poisson manifold in the deformation family constructed above with Poisson structure $\sigma_a$ which is generically symplectic. Then  the cohomology class  in $ H^2(M\backslash D,\C)$ of the dual meromorphic $2$-form $\varphi_a$ is $[\varphi_a]=[\varphi_0]-2a[\omega]$.
\end{prp}

\begin{prf}  
From Proposition \ref{coh} we can use our formulas for $t=0$ at $t=a$, since the cohomology class of $\omega_a$ is the same. In the local basis  from (\ref{sigmaa}) we have $$\sigma_a(\xi_i,\xi_j)=(1-\bar\phi\phi)^{-1}_{ki}\sigma(dz_k,dz_{\ell})(1-\bar\phi\phi)^{-1}_{\ell j}$$ but since $\phi$ is of order $t$ then
$\sigma_t(\xi_i,\xi_j)=\sigma^{ij}+O(t^2)$ and likewise
$\varphi=\sigma_{ij}\xi_i\xi_j+O(t^2)$ where $\sigma_{ij}$ is the inverse of $\sigma^{ij}$. By the same token we also have 
$$\xi_i=dz_i-\phi^i_{\bar k}d\bar z_k+O(t^2)=dz_i-t\sigma^{i\ell}\omega_{\ell\bar k}d\bar z_k+O(t^2).$$
Hence
\begin{eqnarray*}
\varphi&=&\sigma_{ij}dz_idz_j+t(\sigma_{ij}\sigma^{i\ell}\omega_{\ell\bar k}dz_jd\bar z_k-\sigma_{ij}\sigma^{j\ell}\omega_{\ell\bar k}dz_id\bar z_k)+O(t^2)\\
&=&\varphi_0-2t\omega+O(t^2.)
\end{eqnarray*}
Thus the derivative of the cohomology class is $-2[\omega]\in H^2(M,\C)$ restricted to $M\backslash D$.

But from Proposition \ref{coh}  the cohomology class of $\omega_a$ is constant, consequently the variation is linear in $t$.
\end{prf}

\begin{rmks} 

\noindent 1. Note that this linear variation lends a natural role to the parameter $t$ in the construction  analogous, from the Duistermaat-Heckman theorem, to the value of the moment map in an abelian symplectic quotient.

\noindent 2. For a hyperk\"ahler manifold the periods of the holomorphic $2$-form under the twistor deformation  define a conic in a plane in $P(H^2(M,\C))$. Contrast this with the projective line which Proposition \ref{pi} shows occurs in  our Poisson deformation. 

\noindent  3. If $t$ and $\omega$ are real then the imaginary part of $\varphi$ is unchanged. If we had computed the Poisson tensor directly instead of $\varphi$ we would have found that the imaginary part of $\sigma_t$ is  unchanged. This is in fact the situation in bihermitian geometry, where $g([I_+,I_-]X,Y)$ defines the imaginary part of holomorphic Poisson structures in the complex structures $I_+,I_-$ \cite{NJH2}.
\end{rmks}

\section{The two-dimensional case}\label{two}
\subsection{Features}\label{feat}
In two dimensions, as we noted, there is no integrability condition for a Poisson structure and so all we need is an effective anticanonical divisor: $\sigma$ is  simply a holomorphic section of the anticanonical bundle $K^*$. Surfaces which admit such a divisor are either tori or K3 surfaces if $\sigma$ is everywhere non-zero, or certain rational or ruled surfaces \cite{BM}. We are mainly interested in the case where $\sigma$  vanishes on a divisor $D$. From the adjunction formula
$$2(g-1)=KD+D^2=-K^2+K^2=0$$
and so if $\sigma$ vanishes on $D$ in a nondegenerate way then $D$ is an elliptic curve.  

There is a more concrete way to see this: $D$ acquires a holomorphic vector field from $\sigma$ called the {\it modular vector field}. On $D$, the derivative of $\sigma$ is a well-defined section of $\Lambda^2T\otimes T^*$ and contracting on the first factor gives a vector field. In local coordinates  with $\sigma=\sigma^{12}$ the vector field is
$$X=\frac{\partial \sigma}{\partial z_2}\frac{\partial}{\partial z_1}-\frac{\partial \sigma}{\partial z_1}\frac{\partial}{\partial z_2}$$
and is tangential to $D$.
As we noted in the proof of Proposition \ref{pi}, $\sigma_t(\xi_i,\xi_j)=\sigma^{ij}+O(t^2)$ and it follows directly that  the $t$-derivative of the modular vector field at $t=0$ is zero, hence our variation of Poisson structure preserves the vector field on $D$. 

\begin{rmk} In two dimensions   the nondegenerate pairing $T\otimes T\rightarrow \Lambda^2T$ yields an isomorphism  $T\otimes K\cong T^*$  and so for an anticanonical divisor $D$, $T^*\cong T(-D)$,  the sheaf of vector fields vanishing on $D$. The image of the map $\sigma:H^1(S,T^*)\rightarrow H^1(S,T)$ can therefore be thought of as Kodaira-Spencer classes for deformations preserving $D$ and its complex structure.  
The tangential aspect of $D$ is unchanged under deformation but its normal bundle in general does change. We shall see this in a more general context next.
\end{rmk}
\vskip .25cm

 From \cite{BM} if $D$ is nonempty the surface $S$ is ruled or rational and if $H^1(S,{\mathcal O})=0$ it is rational. We restrict now  to rational surfaces.  

Consider the exact sequence of sheaves
$$0\rightarrow K\stackrel{\sigma}\rightarrow{\mathcal O}\rightarrow {\mathcal O}_D\rightarrow 0$$
  Since for a rational surface $H^1(S,{\mathcal O})=0=H^2(S,{\mathcal O})$ we get  from the long exact cohomology 
  sequence
  $$\C\cong H^2(S,K)\cong H^1(D,{\mathcal O}).$$
  From the Dolbeault point of view this isomorphism can be seen as follows. We represent a class in  $H^2(S,K)$ by a $(2,2)$-form $\nu$ and contract with $\sigma$ to get a $\bar\partial$-closed $(0,2)$-form $\sigma(\nu)$. Since $H^2(S,{\mathcal O})=0$ we write this as $\bar\partial\theta$ for a $(0,1)$-form $\theta$. Restricting to $D$, where $\sigma$ vanishes, $\theta$ is $\bar\partial$-closed and represents the class in $H^1(D,{\mathcal O})$.

  Denote this isomorphism by $\alpha: H^2(S,K)\rightarrow H^1(D,{\mathcal O})$. Recall that $D$ and its complex structure are unchanged under our deformation, so $\alpha$ is independent of $t$. 
  
  A holomorphic line bundle $L$ on $S$ is uniquely determined up to isomorphism  by its Chern class $c_1(L)\in H^2(S,\Z)$ since $H^1(S,{\mathcal O})=0$. On the other hand since $H^1(D,{\mathcal O})\cong \C$ its restriction to $D$ has deformations. So we can ask how this restriction varies with $t$ under deformation. 
  
  \begin{prp} \label{var} Let $\omega$ be a closed $(1,1)$-form and consider the deformation given by $\sigma(\omega)$. For each complex structure in the deformation let $L$ be the holomorphic line bundle with first Chern class $c_1(L)$. Then the first variation in the holomorphic structure of  $L$ restricted to $D$ is given by $2\pi i\alpha([\omega] c_1(L))\in H^1(D,{\mathcal O})$, identifying $H^2(S,K)$ with $H^4(S,\C)$.
  \end{prp}
  
  \begin{prf} We again use a \v{C}ech approach. If $g_{UV}$ is a set of holomorphic transition functions for $L$ then 
  $\bar\partial_tg_{UV}=0$. Differentiating at $t=0$ we have 
 \begin{equation}
 \bar\partial\dot g_{UV}+\phi_1(g_{UV})=0
 \label{g}
 \end{equation}
   On $D$, $\sigma=0$ so $\phi_1=0$ and  $\bar\partial\dot g_{UV}=0$. Then $g_{UV}^{-1}\dot g_{UV}$ is a holomorphic cocycle on $D$ representing in $H^1(D,{\mathcal O})$ the first variation of $L$.

  A metric on the holomorphic line bundle $L$ gives rise on $S$ to local functions $h_U$ where $\partial\bar\partial h_U=F\vert_U$ where $F$, the curvature, is a closed $(1,1)$-form. On $U\cap V$ we have 
  \begin{equation}
  \partial h_V-\partial h_U=g_{UV}^{-1}\partial g_{UV}
  \label{h}
  \end{equation}
 
  Now since $\bar\partial\phi_1=0$, 
  $$\bar\partial(\phi_1(h_U))=-\phi_1\bar\partial h_U=-\sigma^{i\ell}\omega_{\ell \bar j}\frac{\partial^2 h_U}{\partial z_i\partial \bar z_k}d\bar z_j d\bar z_k=-\sigma^{i\ell}\omega_{\ell \bar j}F_{i \bar k}d\bar z_j d\bar z_k.$$
 But this is $\bar\partial\theta$ where $\theta$ restricted to $D$ represents $2\pi i\alpha([\omega] c_1(L))$, so $\bar\partial( \phi_1(h_U)-\theta)=0$ and there exists $f_U$ such that
 \begin{equation}
  \phi_1(h_U)-\theta=\bar\partial f_U
  \label{theta}
  \end{equation}
But from (\ref{h}) and (\ref{g})
$$\phi_1(h_V)-\phi_1(h_U)=g_{UV}^{-1}\phi_1(g_{UV})=-g_{UV}^{-1}\bar\partial \dot g_{UV}=-\bar\partial(g_{UV}^{-1}\dot g_{UV})$$
since $g_{UV}$ is holomorphic.
Hence 
$$\bar\partial f_V-\bar\partial f_U=-\bar\partial(g_{UV}^{-1}\dot g_{UV})$$
and so $f_V-f_U+g_{UV}^{-1}\dot g_{UV}$ is a holomorphic $1$-cocycle. But $H^1(S,{\mathcal O})=0$ so there exist local holomorphic functions $a_U$ such that 
$$(f_V-a_V)-(f_U-a_U)+g_{UV}^{-1}\dot g_{UV}=0.$$
Now restrict $\theta$ to $D$. From (\ref{theta}) $\theta=-\bar\partial f_U=-\bar\partial (f_U-a_U)$. Hence its \v{C}ech representative on $U\cap V$ is
$$-(f_V-a_V)+(f_U-a_U)=g_{UV}^{-1}\dot g_{UV}$$
as required.
  \end{prf}
  
  The normal bundle of $D$ is the restriction of $K^*$ hence from Proposition \ref{var}, the normal bundle varies non-trivially if $c_1(K) [\omega]\ne 0$. Consider now some examples.
  
  \subsection{The projective plane}
  
  The projective plane $\PP_2$ is of course rigid so one might question whether the deformation theorem gives any information. But it gave us a deformation of a holomorphic Poisson structure, and this is by no means unique. In the generic case, the anticanonical divisor $D$ is a smooth cubic curve. If $H$  is the hyperplane divisor  then $c_1(K)=-3[H]$ and $[\omega]=k[H]$ so $c_1(K)[\omega]=-3k$ is always nonzero, hence from Proposition \ref{var} the normal bundle of $D$ varies as we vary $t$. We get different embeddings of the same elliptic curve as a plane cubic. 
  
  The viewpoint we get here of a fixed $C^{\infty}$ manifold with different complex structures is close to that in \cite{NJH} for Del Pezzo surfaces. In that paper $[\omega]=c_1$. The deformation of Poisson structure is obtained by putting a metric on $K^*$ and using the function $\log \Vert \sigma \Vert^2$ to define a real Poisson vector field using the real part of $\sigma$. This extends  to a translation on the elliptic curve $D$ and integrating it to a diffeomorphism $f$, the new Poisson structure is $f^*\sigma$ with respect to the transformed complex structure (equivalent of course by $f$ to the original one). 
  
   As remarked in Section \ref{def}, taking  $[\omega]=c_1$ always results  in a trivial Kodaira-Spencer class. Moreover, from Proposition  \ref{coh} the cohomology class of $\omega_t$ is constant, so  all such Kodaira-Spencer classes are trivial and the deformation itself is trivial.   
   
   \subsection{A generic rational surface}
   
   Let $S$ be the surface obtained by blowing up $n$ points $x_1,\dots,x_n$ in $\PP_2$, and let $E_i$ be the divisors of the exceptional curves. Then if $p:S\rightarrow \PP_2$ is the projection, 
   $K_S\sim p^*K+\sum_1^nE_i$ and so 
   $$-K_S\sim p^*3H-\sum_1^nE_i$$
  where $H$ is the hyperplane divisor on $\PP_2$.  A generic effective anticanonical divisor $D$ is therefore the proper transform of a nonsingular cubic curve $C$ passing through the points $x_1,\dots, x_n$. 
   
   Take the class $[\omega]=\sum_1^nm_i[E_i]$ and consider the deformation.  From Proposition \ref{var} the first variation of the line bundle with divisor $E_i$ is
   $\alpha([\omega] [E_i])=m_iu$ where $u\in H^1(D,{\mathcal O})$ is a fixed generator. (In fact, there is a natural one since the modular vector field $X$  gives a trivialization of $K^*_D$ so $H^1(D,{\mathcal O})\cong H^0(D,K_D)^*\cong H^0(D,K^*_D)$.) 
   
   Now curves of self-intersection $-1$ are preserved under deformation so for each $t$ we have a divisor $E_i(t)$ in the same cohomology class which meets $D$ at a point $p_i(t)$, and this point uniquely determines the divisor class on $D$. The deformed complex structure therefore consists of  blowing up points $x_i(t)$ on $C$ moving with uniform velocity $m_iX$.

 When $n>9$, $c_1^2<0$ and so the restriction of $K^*$ to  $D$ has negative degree. It follows that $\dim H^0(S,K^*)=1$ and the surface has (up to a multiple) a unique Poisson structure.  Also, 
   in two dimensions  $\sigma$ is a section of the line bundle $K^*$ so the map  $\sigma:H^1(S,T^*)\rightarrow H^1(S,T)$ appears in  the cohomology sequence of the exact sequence of sheaves 
\begin{equation}
0\rightarrow T^*\stackrel{\sigma}\rightarrow T\rightarrow T\vert_D\rightarrow 0
\label{Texact}
\end{equation}
that is,
$0\rightarrow H^0(S,T)\rightarrow H^0(D,T)\rightarrow H^1(S,T^*)\stackrel{\sigma}\rightarrow H^1(S,T)\rightarrow\cdots$
It follows that the kernel of $\sigma$ is  isomorphic to $H^0(D,T)$  since if $n>4$ then $H^0(S,T)=0$. The normal bundle $K^*$ of $D$ has negative degree so all sections of $T$ on $D$ are tangential to $D$ and hence form a one-dimensional space. In this case, then, the kernel of $\sigma:H^1(S,T^*)\rightarrow H^1(S,T)$ is one-dimensional. We know   $c_1$  is always in the kernel, so it is the generator. The image of $\sigma$ is thus the image of $c_1^{\perp}$ relative to the intersection form. Note that   this has the integral structure of the Dynkin diagram $E_n$.

From Proposition \ref{pi} the periods of the unique meromorphic one-form determine, at least locally,  the modulus of such Poisson surfaces.

   \subsection{A special rational surface}
   In the previous example, the subspace of $H^1(S,T^*)$ which gave trivial deformations  was one-dimensional, the smallest possible. The next example is the opposite extreme.
   
   We produce a Poisson surface by  blowing up points in a highly specialized fashion, following \cite{Sak}. 
   Take a nonsingular cubic curve $C$  in $\PP_2$ and a  line $L$ tangent to $C$ at an inflection point $x\in C$. Now proceed to blow up three times (see \cite{Sak} p.222 for explicit formulas) taking as centre  each time the  point of intersection of the proper transform of $L$ with the exceptional curve. At this stage the second order tangency of $C$ with $L$ yields another distinguished point on the exceptional curve and we blow that up $5$ more times. Finally choose a point to blow up on the exceptional curve.  We get a surface $S$ with $c_1^2=0$ and $c_2=12$ and the collection of $-2$ curves formed in the blowing up process gives
   a configuration of rational curves intersecting according to the extended Dynkin diagram of $E_8$:
   $$   \eddeiii 2.3.4.6.5.4.3.2.1.$$
 There is an anticanonical divisor $D$ supported on these curves where the multiplicity of each component of $D$ is the number on the corresponding node. 
The surface $S$ has the property that each irreducible component $D_i$ of $D$ satisfies   
 $KD_i=0$. 
 
  As remarked in Section \ref{def}, each component  of an anticanonical divisor defines a cohomology class in the kernel of $\sigma:H^1(S,T^*)\rightarrow H^1(S,T)$ so we have at least a $9$-dimensional kernel. On the other hand, blowing up at the fourth stage and beyond kills any holomorphic vector field   and so the final choice of a point $p$ to blow up gives an effective one-parameter family of deformations, preserving the divisor.  Since $c_2=12$, $\dim H^1(S,T^*)=10$ hence the image in $H^1(S,T)$ is one-dimensional.
 
 Now let $E$ be  the exceptional curve  created in the last blow-up.  It has self-intersection $-1$ so $K E=-2-E E=-1$ and  the divisor class $[E]$ is not a linear combination of the $[D_i]$. Hence the only possible non-trivial deformation comes from taking $[\omega]$ to be a multiple of $[E]\in H^1(S,T^*)$. 
 
 Let $S$ be a surface of this type where the point $p$ is not the intersection with the proper transform of $C$.  Then the exceptional curve $E$ meets the last component $D_9$ in a single point which is not the intersection with $D_8$. Since the multiplicity of $D_9$ (the last node in the Dynkin diagram) is $1$, this intersection point determines the divisor class of $E$ restricted to $D$. Similarly, the modular vector field of $\sigma$ is a  vector field on the rational curve $D_9$ with a double zero at $D_8\cap D_9$. Since $K E=-1$ it follows from Proposition \ref{var} that  the deformation $\sigma([E])$ is obtained by moving the final point $p$ along the curve whose proper transform is $D_9$. Furthermore, as in the previous  example the point moves on $\C=D_9\backslash D_8$ with constant speed with respect to the parameter $t$ in the deformation. 
 
  \begin{rmk} This family  of Poisson surfaces arises in \cite{Sak} in the context of the first Painlev\'e equation $y''=6y^2+x$. The parameter $x$ in this equation is essentially $t^{-1}$ in terms of the deformation parameter $t$. The compact surface $S$ undergoes a non-trivial deformation but $S\backslash D$ does not -- in fact the Painlev\'e equation is a time-dependent vector field which integrates to a family $f_t$ of symplectic diffeomorphisms of $S\backslash D$.
   \end{rmk}
   
    \subsection{ALE spaces}

   An ALE space is a non-compact hyperk\"ahler 4-manifold $M$ which is asymptotic  to $\C^2/\Gamma$ with its Euclidean metric, where $\Gamma\subset SU(2)$ is a finite subgroup. They can all be constructed as finite-dimensional hyperk\"ahler quotients \cite{K1}. Infinity is modelled on $\R^4/\Gamma$ and the manifold admits an orbifold conformal compactification. Using the  twistor space $Z$, the   space $M$ with one of its complex structures  compactifies to a singular surface $S$ (see \cite{K2}), which is a Poisson surface. We shall look here at the relationship between deformations of $S$ and the well-understood moduli of ALE hyperk\"ahler metrics. 
      
 The simplest example is obtained from the quotient of $\PP_2=\C^2\cup \PP_1$ by the  extended action of  $\Gamma$ on  $\C^2\oplus \C$. The singularity at $[0,0,1]$ is the origin in $\C^2$ which is resolved to give an ALE space. The fixed points on the line at infinity $[z_1,z_2,0]$ give simple Hirzebruch-Jung singularities on the quotient and  resolving them gives a smooth surface $S$ which is a compactification of $M$ by an anticanonical divisor -- a section of $K^*$ which extends the inverse of the holomorphic symplectic $2$-form $\varphi$ on $M$ which is part of the hyperk\"ahler picture. 
    
    This 
   divisor is  another configuration of rational curves. The simplest ones from this point of view are  given by the binary tetrahedral, octahedral and icosahedral groups where there are three singular  points corresponding to the stabilizers of vertices, edges and faces. The anticanonical divisor  $D$ is described by  the graphs
  $$ \dddiv2.1.3.3. \dddiv2.1.3.4. \dddiv2.1.3.5.$$
with $-K\sim 2C_0+C_1+C_2+C_3$ where $C_0$ is the central curve, and where now the  number $m$ at a vertex corresponds to a rational curve of self-intersection $-m$.

The general ALE space replaces the resolution of the quotient $\C^2/\Gamma$ by its versal deformation but the divisor at infinity is the same. Whereas the simplest case has a $\C^*$ symmetry $[z_1,z_2,z_3]\mapsto  [ z_1,z_2,\lambda z_3]$ this does not hold in general. 

There is a Torelli theorem for ALE spaces \cite{K2} which describes explicitly the parameters for deforming the hyperk\"ahler metric. These describe a deformation of the compactification $S$, and if we are only interested in the complex  structure it is  the periods of the $2$-form $\varphi$ which determine it --  its cohomology class in $H^2(M,\C)$, which  has the structure of the Cartan subalgebra of type $E_6,E_7$ or $E_8$. The Euler characteristic of $D$ is $5$ so $\dim H^2(S,\C)=\dim H^1(S,T^*)=10,11$ or $12$. 

If we now use the exact sequence (\ref{Texact}) we have 
   $0\rightarrow H^0(S,T)\rightarrow H^0(D,T)\rightarrow H^1(S,T^*)\stackrel{\sigma}\rightarrow H^1(S,T)\rightarrow\cdots$
   and for a generic ALE,  $H^0(S,T)=0$. We know that the cohomology classes of the four components of $D$ lie in the kernel of $H^1(S,T^*)\stackrel{\sigma}\rightarrow H^1(S,T)$, so if $\dim H^0(D,T)=4$ we deduce that the ALE deformations give effective deformations of the compactification $S$. 
      
  To see that this is true, note that  all components of $D$ have negative normal bundles, and so any section of $T$ on $D$ is tangential to the component. Moreover  since $C_0$ has multiplicity $2$ any section on $C_1,C_2$ or $C_3$ is a vector field on a rational curve which vanishes with multiplicity $2$ at its point of intersection with $C_0$. This provides a one-dimensional space for each $C_i$, $i\ne 0$. On $C_0$ itself we need a vector field vanishing at three points, which must be zero. But $C_0$ has multiplicity $2$ and normal bundle ${\mathcal O}(-1)$ so since $T_{C_0}(1)\cong {\mathcal O}(3)$, there is a one-dimensional space of such sections, giving four dimensions in all. 

\section{Hilbert schemes}\label{HS}

\subsection{Deformations}
The Hilbert scheme $S^{[n]}$ of a surface $S$ parametrizes zero-dimensional subschemes of length $n$ on $S$ and is a resolution of the singularities of the symmetric product $S^{(n)}=S^n/\Sigma_n$ where $\Sigma_n$ is the symmetric group. It inherits many properties  of $S$: in particular if $H^p(S,{\mathcal O})=0$ for $p>0$ then the same is true for $S^{[n]}$.

A Poisson structure $\sigma$ on $S$ defines canonically one on the product  $S^n$
$$\sigma(x_1,x_2)\frac{\partial}{\partial x_1}\wedge\frac{\partial}{\partial x_2}+\sigma(y_1,y_2)\frac{\partial}{\partial y_1}\wedge\frac{\partial}{\partial y_2}+\cdots $$ which is invariant by the symmetric group. It is shown in \cite{Bot} that this extends canonically to a Poisson structure $\tau$ on  the Hilbert scheme. This structure is generically symplectic and if  $\sigma$ vanishes on a  smooth elliptic curve $D$ then $\tau$ vanishes on the smooth symmetric product $D^{(n)} (=D^{[n]})\subset S^{[n]}$.

If the complex structure of $S$ is deformed then there is a corresponding deformation of $S^{[n]}$ and the more general question of the relation between the two deformation functors was addressed by Fantechi in \cite{Fan}. We briefly describe the approach (referring to \cite{Fan} for more details).

Let $p:S^{[n]}\rightarrow S^{(n)}$ be the natural map resolving the singularities then, with $\theta_X$ denoting the tangent sheaf,  $p_*\theta_{S^{[n]}}=\theta_{S^{(n)}}$ and the Leray spectral sequence gives
$$0\rightarrow H^1(S^{(n)},\theta)\rightarrow H^1(S^{[n]},\theta)\stackrel{\rho}\rightarrow H^0(S^{(n)},R^1p_*\theta)\rightarrow H^2(S^{(n)},\theta)\rightarrow \cdots $$
The term $H^1(S^{(n)},\theta)$ is the invariant part under the $\Sigma_n$-action of $H^1(S^n,\theta)$ and if $H^1(S,{\mathcal O})=0$ then by the K\"unneth formula this is  naturally isomorphic to $H^1(S,\theta)$. The first  part of  the sequence  can then be written as $0\rightarrow H^1(S,\theta)\rightarrow H^1(S^{[n]},\theta)$ and is the natural map of Kodaira-Spencer classes  for an infinitesimal deformation of $S^{[n]}$ induced by one of $S$. 

Now let $S_{\mathrm {sing}}^{(n)}$ denote the singular locus of $S^{(n)}$. Its resolution in $S^{[n]}$ is the {\it exceptional divisor} $E$.  
\begin{rmk}
The exceptional divisor  defines a distinguished class $[E]\in H^2(S^{[n]},\Z)$. Together with the classes pulled back from $S^{(n)}$ (which correspond to those in $S$) $[E]/2$  generates the second cohomology. 
\end{rmk}
Let $Z$ be the  smooth locus of $S_{\mathrm {sing}}^{(n)}$, where just two points coincide.  Then $S_{\mathrm {sing}}^{(n)}$ has codimension $2$ in $S^{(n)}$ and the singular set of  $S_{ \mathrm{sing}}^{(n)}$ is of codimension $2$ in  $S_{\mathrm {sing}}^{(n)}$. The two coincident points define a projection $q: Z\rightarrow S$ and on $Z$ one can show that $R^1p_*\theta\cong q^*(K^*_S)$. Taking account of the codimensions, Hartogs' theorem gives 
$$H^0(S^{(n)},R^1p_*\theta)\cong H^0(S,K^*).$$
The appearance of the space of Poisson tensors on $S$ in the computation of $H^1(S^{[n]},\theta)$  is suggestive, and explained by the theorem below, which shows in particular that we have a split exact sequence 
$$0\rightarrow H^1(S,T)\rightarrow H^1(S^{[n]},T)\rightarrow H^0(S,K^*)\rightarrow 0$$
(since $S^{[n]}$ is smooth we go back to using $T$ instead of $\theta$).

\begin{thm} \label{hilbdef} Let $S$ be a surface with $H^1(S,{\mathcal O})=0$ and let $[E]\in H^1(S^{[n]}, T^*)$ be the cohomology class of the exceptional divisor on the Hilbert scheme $S^{[n]}$. Let $\sigma$ be a Poisson structure on $S$ and $\tau$ the induced one on $S^{[n]}$. Let  $\rho: H^1(S^{[n]},T)\rightarrow H^0(S,K^*)$ be the homomorphism above, then 
$\rho\tau([E])=-2\sigma.$
\end{thm}

\begin{prf} We first need to consider in more detail the isomorphism in \cite{Fan} between $R^1p_*\theta$ on $Z$ and $q^*(K^*_S)$. In a neighbourhood of a point in $Z$ only two points coalesce and so the resolution locally looks like a product of an open set in $\C^{2n-4}$ and $S^{[2]}$. The Hilbert scheme $S^{[2]}$ has a concrete construction -- blow up the diagonal in $S^2$ and divide out by the involution interchanging the two factors. Since the fixed point set is the exceptional divisor of the blow-up the quotient is smooth. 

The normal bundle of the diagonal $\Delta\subset S^2$ is the tangent bundle $T_S$ and so the blow-up replaces $\Delta$ by the projective bundle $\pi: \PP(T)\rightarrow S$ and its normal bundle $L$ is the tautological bundle $L\subset  \pi^*T$. Taking the quotient by the involution gives normal bundle $N=L^2$. 

Now local sections of  $R^1p_*\theta$ are non-zero only on a neighbourhood of a point on the diagonal so consider such a neighbourhood $U$. Then $H^1(p^{-1}(U),T)$  is defined as sections  over $U\cap Z$ of  the vector bundle whose fibre over $x\in S\cong  \Delta$ is $H^1(\PP(T_x),T)$. But the tangent bundle of $\PP(T)$ restricted to a fibre is ${\mathcal O}(2)\oplus {\mathcal O}^2$ hence $H^1(\PP(T_x),T)\cong H^1(\PP(T_x),N)=H^1(\PP(T_x),L^2)$. But from the Euler sequence $K^*_x=\Lambda^2T_x$ is naturally identified with $K_F^*L^2$ where $K_F$ is the canonical bundle along the fibres of $\pi$, hence 
$$H^1(\PP(T_x),L^2)\cong K_x^*\otimes H^1(\PP(T_x),K_F)\cong K_x^*$$
where the last isomorphism comes from taking the standard trivialization of $R^1\pi_*K_F$. 

Now consider the behaviour of the Poisson tensor in the resolution.  In local coordinates the relevant piece on $S\times S$ is 
\begin{equation}
\sigma(x_1,x_2)\frac{\partial}{\partial x_1}\wedge\frac{\partial}{\partial x_2}+\sigma(y_1,y_2)\frac{\partial}{\partial y_1}\wedge\frac{\partial}{\partial y_2}.
\label{sig}
\end{equation}
Writing $u_i=x_i-y_i, v_i=x_i+y_i$ the involution exchanging factors is $(u,v)\mapsto(-u,v)$ and the Hilbert scheme is given by resolving the singular quotient  $\C^2/\pm 1$. Writing $x=u_1^2,y=u_2^2,z=u_1u_2$ the singularity is the cone $xy=z^2$ which is resolved by the total space of the line bundle ${\mathcal O}(-2)$ on $\PP_1$ (the cotangent bundle): the singular origin in $\C^3$ is replaced by the zero section of ${\mathcal O}(-2)$.  If $\zeta$ is an affine coordinate on $\PP_1$ and $\eta d\zeta$ the cotangent vector at a point on $T^*\PP_1$ then in these local coordinates $\eta=u_1^2,\zeta=u_2/u_1$ and 
$$
2\frac{\partial}{\partial u_1}\wedge\frac{\partial}{\partial u_2}=\frac{\partial}{\partial \eta}\wedge\frac{\partial}{\partial \zeta}
$$
which is the standard Poisson structure of the canonical symplectic structure on $T^*\PP_1$. 

Hence if $f(u,v)$ is invariant by the involution $(u,v)\mapsto(-u,v)$, then $f{\partial}/{\partial u_1}\wedge{\partial}/{\partial u_2}$ extends on the resolution, and the Poisson tensor (\ref{sig}) can be written as
$$f \frac{\partial}{\partial \eta}\wedge\frac{\partial}{\partial \zeta}+\frac{\partial}{\partial \eta}\wedge\left(a_1\frac{\partial}{\partial v_1}+a_2\frac{\partial}{\partial v_2}\right)+\frac{\partial}{\partial \zeta}\wedge\left(b_1\frac{\partial}{\partial v_1}+b_2\frac{\partial}{\partial v_2}\right)+c\frac{\partial}{\partial v_1}\wedge\frac{\partial}{\partial v_2}.$$

We want to calculate $\rho\tau([E])$. This means first taking $[E]\in H^1(S^{[n]}, T^*)$ and restricting to $\PP(T_x)$, which is $\eta=0$ in our local coordinates. This lies in the one-dimensional space $H^1(\PP(T_x),K_F)$ and is $-2 \times $ the standard generator. A \v Cech representative is of the form $a(\zeta)  d\zeta$.
Now take the Poisson tensor $\tau$ and contract to get a class in $H^1(\PP(T_x),T)$. This is one-dimensional and is isomorphic to $H^1(\PP(T_x),N)$. A \v Cech representative for this is of the form $b(\zeta) \partial/\partial\eta$.  It follows that in the above expression for a Poisson tensor on the resolution  only the   first term gives a contribution. The coefficient of this is $f$, and changing to the original coordinates this is   $\sigma(x_1,x_2)$. 
\end{prf}

This theorem  tells us in particular that $\sigma([E])\in H^1(S^{[n]}, T)$ is non-zero, so applying the deformation results of Section \ref{def} we have a one-parameter family determined by a section of $K^*$ on $S$. 

\begin{ex} Take $S=\PP_2$ and the divisor $D$ to be $3L$ where $L$ is the line at infinity. The deformations of the Hilbert scheme preserve the open subset on which $\tau$ is symplectic and thus give  deformations of the Hilbert scheme of $\C^2$. These have explicit descriptions as finite-dimensional hyperk\"ahler quotients \cite{Nek}. 
\end{ex}

In Section \ref{two} we examined the variation in the holomorphic structure of a line bundle in a given cohomology class under the deformation. We do this next for two particular classes in $H^2(S^{[n]},\Z)$ -- the first Chern class $c_1$ of the manifold and the class $[E]$ of the exceptional divisor.

\begin{prp} \label{line} Let $(S,\sigma)$ be a rational Poisson surface such that $\sigma$ vanishes on a smooth elliptic curve $D$. Then, in the deformation of the Hilbert scheme $S^{[n]}$ produced in Theorem \ref{deform} by taking $[\omega]=[E]$, the class of the canonical bundle restricted to the zero set $D^{(n)}$ of the Poisson structure  is constant  in $H^1(D^{(n)},{\mathcal O}^*)$, but the line bundle with Chern class $[E]$ varies linearly in $t$.
\end{prp}
\begin{prf} We shall reduce the question to the case of a surface and then use the method of Proposition \ref{var}. 

We start with the Hilbert scheme itself. A point $a$ lying on the curve $D\subset S$ defines $a^{n-1}\in D^{(n-1)}\subset S^{[n-1]}$ corresponding in coordinates to the ideal $\{z_1,z_2^{n-1}\}$ where $a$ is the origin and $z_1=0$ a local equation for $D$. We can consider  in $S^{[n]}$ the subschemes of length $n$ containing this. 
Clearly for $x\in S$, $x\ne a$ there is such a scheme supported on $x$ and $a$ and this defines an embedding of $S\backslash\{a\}\subset S^{[n]}$ which extends as $x$ approaches $a$ to the blow-up $\hat S$ of $S$ at $a$. The blow up consists of ideals $\{z_1,(\lambda_1z_1+\lambda_2z_2)z_2^{n-1}\}$ supported at $a$. The map $p:S^{[n]}\rightarrow S^{(n)}$ blows down the exceptional curve on $\hat S$ and the image is  $S\cong S\times \{a\}\times  \dots\times \{a\}\subset S^{(n)}$.

This holds for any curve in $S$, but when we take $D$, then the induced Poisson structure $\tau$ is tangential to $\hat S$ since $\sigma$ vanishes at $a$. In other words, there is a Poisson map from $\hat S$ (with the  Poisson structure determined  by the proper transform $\hat D$ of the cubic curve $D$ through $a$) to $S^{[n]}$. The deformation of Theorem \ref{deform} therefore preserves $\hat S$, though not necessarily its complex structure -- it is a deformation as in Theorem \ref{deform} for the surface $\hat S$ with induced Poisson structure $\hat\sigma$ vanishing on $\hat D$.

The structure of the zero set of $\tau$, namely the symmetric product $D^{(n)}$, can be seen by associating to an $n$-tuple of points its divisor class. This  represents $D^{(n)}$ as a projective bundle over $\Pic^n(D)$ and so $H^1(D^{(n)},{\mathcal O})\cong H^1(\Pic^n(D),{\mathcal O})$. The surface  $\hat S$ intersects $D^{(n)}$ in $\hat D$. The map $x\mapsto [x+(n-1)a]$ identifies $D$ with $\Pic^n(D)$  so we have natural identifications of $H^1(D^{(n)},  {\mathcal O})$ and $H^1( D,  {\mathcal O})=H^1(\hat D,  {\mathcal O})$. Using Proposition \ref{var} for the surface $\hat S$ it suffices to prove that, if $E$ is the restriction of the exceptional divisor  on $S^{[n]}$ to $\hat S$, then $c_1(E)c_1(K_{S^{[n]}}) =0$ and  $c_1(E)^2\ne 0$. 

 The Chern class of the exceptional divisor $E$ can be determined from Lemma 3.7 in \cite{Lehn}. There is a universal sheaf $\Xi$ on $S^{[n]}\times S$. If  $\pi$ is the projection onto the first factor, then it is finite of degree $n$  so $\pi_*({\mathcal O}_{\Xi})$ is a rank $n$ vector bundle. The formula then is 
 $$c_1(E)=-2c_1(\pi_*{\mathcal O}_{\Xi}).$$
 If $C$ is a curve in $S$ which does not meet $a$ then it lifts to $\hat S$, and  the vector bundle  restricted to $C$ is the direct sum of a trivial rank $(n-1)$ bundle corresponding to $a^n$ and the direct image of the diagonal $\Delta$ in $C\times C$.   From Grothendieck-Riemann-Roch   
 \begin{equation}
 c_1(E)[C]=4g-4.
 \label{EC}
 \end{equation}
 The canonical bundle of the Hilbert scheme  $S^{[n]}$  is the pull-back $p^*K_{S^{(n)}}$ from the symmetric product and so restricted to $\hat S$ this is $p^*K_S$. Since a divisor of $K^*_S$  is an elliptic curve,  $g=1$ and we have  $c_1(E)c_1(K_{S^{[n]}}) =0$.
 
 The exceptional curve $F$ in $\hat S$ obtained by blowing up $a$ is the projective line with homogeneous coordinates $[\lambda_1,\lambda_2]$ where the ideal is $\{z_1,(\lambda_1z_1+\lambda_2z_2)z_2^{n-1}\}$.  It follows that on $\PP_1\times \PP_1$  the sheaf is the divisor $\Delta +(n-1)\{[0,1]\times \PP_1\}$ which is in the divisor class ${\mathcal O}(1,n)$.  By Grothendieck-Riemann-Roch again we have 
  \begin{equation}
 c_1(E)[F]=-2(n+1).
  \label{EF}
 \end{equation}
 To calculate $c_1(E)^2$ it is sufficient to take $S$ to be $\PP_2$ blown up at $(k-1)$ distinct points,  in which case, applying (\ref {EC}) and (\ref{EF}) to a line and the exceptional curves 
 $$c_1(E)^2=32-16k-4(n+1)^2$$
 which is always negative for $n>1$. 
 \end{prf}

This proposition identifies a geometric object which is changing  under the deformation. Whereas the holomorphic structure of the zero set of the Poisson tensor is unchanged, the cohomology class of the exceptional divisor on that subvariety varies with $t$.

\subsection{The Hilbert scheme $\PP_2^{[n]}$}

Since $\PP_2$ is rigid, Theorem \ref{hilbdef} gives
$$ H^1(\PP_2^{[n]},T)\cong H^0(\PP_2,K^*)$$
and Theorem \ref{def} tells us that the Hilbert scheme $\PP_2^{[n]}$ has non-trivial deformations, each one a Poisson manifold. Consider the generic case where the initial Poisson structure  is induced by a section of $K^*$ on $\PP_2$ which vanishes on a  smooth  cubic curve $D$. The zero set of the Poisson tensor on $\PP_2^{[n]}$ is the symmetric product  $D^{(n)}$ and, as we have seen, this is preserved  under our deformation.  

\begin{prp} \label{h1} Let $\sigma$ be a Poisson structure on $\PP_2$ whose zero set  is a  smooth cubic curve, and let $M$ be the deformation of the Hilbert scheme for $t\ne 0$. Then  

\noindent (i) $H^0(M,T)=0$ 

\noindent (ii) $\dim H^1(M,T)=2$. 
\end{prp}
\begin{prf} 

\noindent (i) First note  that since $p_*\theta_{S^{[n]}}=\theta_{S^{(n)}}$, 
 holomorphic vector fields on the Hilbert scheme  are defined by $\Sigma_n$-invariant vector fields on $S^n=\PP_2^n$. Hence they are all given by the induced action from $\PP_2$. 
 We shall show that there is a first order obstruction to extending any such vector field.
 
 Let $Z$ be a holomorphic vector field to first order in the deformation, with $Z(0)=X$. Then the $T^{1,0}$ component in the  complex structure at $t=0$ is of the form $X+tY$ where
 $$\bar\partial Y+[\tau(\omega),X]=0$$
 which implies that the cohomology class ${\mathcal L}_X(\tau([\omega]))$ is zero in $H^1(\PP_2^{[n]},T)$. 
But ${\mathcal L}_X$ acts trivially on $H^1(\PP_2^{[n]},T^*)$, so this class is $({\mathcal L}_X\tau)([\omega])$. 

However, from Theorem \ref{hilbdef}, this class is zero if and only if  ${\mathcal L}_X\tau=0$. But a smooth cubic is not fixed by any projective transformation, so the class is non-zero and the deformation is obstructed, whatever the choice of $X$.
 
\noindent (ii) Similarly consider a first order extension  of a class in $H^1(T)$. From Theorem \ref{hilbdef} we can represent this by $\pi(\omega)$ on $\PP_2^{[n]}$ for some holomorphic Poisson tensor  $\pi$ induced from an anticanonical section on $\PP_2$, and a first order deformation defines $\alpha \in \Omega^{0,1}(T)$   such that  $\bar\partial_t(\pi(\omega)+t\alpha)=0+O(t^2)$ and hence 
\begin{equation}
\bar\partial\alpha+[\tau(\omega),\pi(\omega)]=0.
\label{obs}
\end{equation}
Here $[\tau(\omega),\pi(\omega)]$ represents the  obstruction class in $H^2(\PP_2^{[n]},T)$ to making the extension, but in fact this vanishes. Recall in the proof of Theorem  \ref{def} (and with $H^2(M,{\mathcal O})=0$) we had 
$$[\sigma(\omega),\sigma(\omega)]=-2\bar\partial (\sigma\partial\beta_2)$$
where $-2\bar\partial \beta_2=\sigma(\omega^2)$. 
But the induced Poisson structures on the Hilbert scheme are linear in the sections of $H^0(\PP_2,K^*)$, so applying this to $\tau,\pi$ and $\tau+\pi$ gives an $\alpha$ satisfying Equation \ref{obs}. 

A cohomology class to the first order is defined by  such a form  modulo $\bar\partial_t(X+tY)$ mod $t^2$, so $\pi(\omega)+t\alpha_1,\pi(\omega)+t\alpha_2$ satisfying (\ref{obs}) define the same class if 
$$\alpha_1-\alpha_2=\bar\partial Y+[\tau(\omega),X]\quad {\mathrm {and}} \quad \bar\partial X=0,$$
or equivalently $[\alpha_2-\alpha_1]=[({\mathcal L}_X\tau)(\omega)]\in H^1(\PP_2^{[n]},T)$.

 But this equivalence is the  $10$-dimensional space  $H^1(\PP_2^{[n]},T)$ modulo the  Lie derivative action of the $8$-dimensional Lie algebra $H^0(\PP_2^{[n]},T)$ on $\tau(\omega)$. Since $\tau$ is not fixed by any vector field, the quotient is $2$-dimensional. Hence $H^1(M,T)$ for a generic deformation is at most $2$-dimensional. 

On the other hand we have two clear parameters in the deformations -- the modulus of the cubic curve which gives the Kodaira-Spencer class at $t=0$,  and, from Proposition \ref{line} the  class of the line bundle determined by the exceptional divisor on $D\subset D^{(n)}$. Since the complex structure of the zero set $D^{(n)}$ is unchanged in the deformation, so is the  modulus of $D$, and hence these  two parameters are independent. We deduce that $\dim H^1(M,T)=2$ generically. 
\end{prf}

 This proposition tells us that a generic deformation has a two-dimensional local moduli space and we have  identified  two parameters. In the case $n=2$ these parameters become more explicit.

\subsection{The case $\PP_2^{[2]}$}

 A pair of unordered distinct points in $\PP_2$ defines a line and hence a point of the dual projective space $\PP^*_2$. Moreover when two points coincide the Hilbert scheme captures the direction (this is what blowing up the diagonal does in the earlier description of $S^{[2]}$) so we get a well defined projection $p: \PP_2^{[2]}\rightarrow \PP_2^*$. The fibre is the symmetric product 
$\PP_1^{(2)}$ of the line in $\PP_2$ dual to $p(x)\in \PP_2^*$. It follows that $\PP_2^{[2]}$ is the projective bundle $\PP(\Sym^2T)$ over $\PP_2^*$. We can also write $\Sym^2T(-3)=\End_0 T$ where $\End_0 $ is the sheaf of trace zero endomorphisms. From this point of view, the two eigenspaces of $A\in \End_0T$ determine two lines in $\PP_2^*$ through $p(x)$, or dually the two points in $\PP_2$ lying on the line dual to $p(x)$. For convenience we shall write $E_0=\End_0 T$. 

The rank $3$ vector bundle $E_0$ has $c_1(E_0)=0,c_2(E_0)=3$ and is stable since  the Fubini-Study metric is K\"ahler-Einstein.

\begin{prp} \label{def2} $H^1(\PP_2^{[2]},T)\cong H^1(\PP^*_2,\End_0 E)$
\end{prp}
\begin{prf} Consider $\PP_2^{[2]}\cong \PP(E_0)$ and view the tangent bundle as an extension
\begin{equation}
0\rightarrow T_F\rightarrow T\rightarrow p^*T_{\PP_2^*}\rightarrow 0
\label{Tseq}
\end{equation}
and use the Leray spectral sequence of the fibration. We obtain  directly that the spaces $H^k(\PP(E_0),p^*T_{\PP_2^*})$ and $H^k(\PP_2^*,T)$ are isomorphic and zero unless $k=0$ and similarly $H^k(\PP(E_0),T_F)\cong H^k(\PP_2^*,p_*T_F)$. But  the fibrewise Euler sequence
$$0\rightarrow {\mathcal O}\rightarrow L^*\otimes E_0\rightarrow T_F\rightarrow 0$$
(where $L$ is the tautological bundle) gives 
$p_*T_F\cong \End_0 E_0$, and hence  $H^1(\PP(E_0),T_F)\cong H^1(\PP_2^*,\End_0 E_0)$. The result follows from the long exact cohomology sequence of (\ref{Tseq}):  by stability, $H^0(\PP_2,\End_0 E_0)=0$ and we have 
\begin{equation}
 0\rightarrow  H^0(\PP(E_0),T)\rightarrow H^0(\PP_2^*,T)\rightarrow H^1(\PP_2^*,\End_0 E_0)\rightarrow H^1(\PP(E_0),T)\rightarrow 0
 \label{ex2}
 \end{equation}
 But  $E_0=\End_0 T$ is acted on naturally by any automorphism of $\PP_2^*$ so the map $H^0(\PP(E_0),T)\rightarrow H^0(\PP_2^*,T)$ is surjective, hence the result.  
\end{prf}

By Riemann-Roch the stable bundle $E_0$ has a smooth moduli space of dimension $\dim H^1(\PP^*_2,\End_0E_0)=10$.   It follows from the Proposition that  our deformations of the Hilbert scheme are all projective bundles $\PP(E)$ over $\PP_2^*$. 

\begin{rmk} The exact sequence (\ref{ex2}) above holds for any of these vector bundles $E$. As in the proof  of Proposition \ref{h1}, if we deform in the direction determined by a non-singular cubic curve, the holomorphic  vector fields on $\PP_2^*$ do not lift to $\PP(E)$ and then it follows directly from the sequence that $\dim H^1(\PP(E),T)=2$.
\end{rmk}

Let ${\mathcal O}(1)$ be the hyperplane bundle on $\PP_2^*$ and $H$ the dual of the tautological line bundle on $\PP(E)$, then $K_{\PP(E)}^*\cong H^3(3)$. The exceptional divisor on $\PP(\End_0T)$, where  two points in $\PP_2$ coincide, is where the two eigenspaces of $A\in \End_0T$ coincide, i.e. where $\tr A^2=0$. It follows that  the exceptional divisor is defined by a section of $H^2$. 

There is a classical  description of the symmetric product $\PP_2^{(2)}$. A pair of distinct points in $\PP_2$ determine  dually a pair of lines in $\PP_2^*$, which is a singular conic. The symmetric product can then be identified with the cubic fourfold $\det S=0$ in the $5$-dimensional projective space of symmetric $3\times 3$ matrices. This has a singularity along the rank one symmetric matrices (where the conic is a double line) and its resolution is the exceptional divisor.

 The projective embedding of a cubic fourfold is given by $K^{-1/3}$, so since  $K^*\cong H^3(3)$ the map from the  Hilbert scheme to the symmetric product is  the map to $\PP_5$ given by sections of the line bundle $H(1)$ on $\PP(\End_0 T)$, or equivalently  sections of $(\End_0T)^*(1)\cong \End_0T (1)$ on $\PP_2^*$. Now $H^p( \PP_2^*, \End_0T(1))=0$ for $p=1,2$ and has dimension $6$ for $p=0$. Vanishing will hold for small deformations $E$ of $E_0=\End_0T$ and hence the line bundle $H(1)$ for deformations of the Hilbert scheme will map $\PP(E)$ to $\PP_5$ and its image will be another cubic fourfold. In terms of the vector  bundle $E$, sections of  $H(1)$ on $\PP(E)$ are naturally isomorphic to sections of $E^*(1)$ on $\PP_2^*$.  
 
   By stability  and $c_1(E)=0$ we have $H^0(\PP_2^*,E)=H^0(\PP_2^*,E^*)=0$ hence in particular $H^0(\PP_2^*,E^*(-1))=0$. It follows from  (\cite{OSS} page 252)  that $E$  is the cohomology of a monad 
$$0\rightarrow H^1(\PP_2^*,E(-2))\otimes {\mathcal O}(-1)\rightarrow  H^1(\PP_2^*,E(-1))\otimes \Omega^1(1)\rightarrow H^1(\PP_2^*,E)\otimes {\mathcal O}\rightarrow 0.$$
Since $E^*$ is also a deformation, and this is what we need for $H^0(\PP(E),H(1))$, we work with $E^*$. 

By Riemann-Roch $H^1(\PP_2^*,E^*)=0$ and $  H^1(\PP_2^*,E^*(-2))$ and $H^1(\PP_2^*,E^*(-1))$
both have dimension $3$. Hence $E^*$, and any deformation of it, appears naturally as a quotient
$$0\rightarrow \C^3\otimes {\mathcal O}(-1)\rightarrow \C^3\otimes \Omega^1(1)\rightarrow E^*\rightarrow 0.$$
Thus $E^*$ is defined by a $3\times 3$ matrix with entries in $H^0(\PP_2^*,\Omega^1(2))$. More concretely, if $[x_1,x_2,x_3]$ are homogeneous coordinates for $\PP_2^*$ then a basis for the  global sections of $T^*(2)$ is given by  $\alpha_1=x_2dx_3-x_3dx_2$ etc. Let $A_{ijk}\alpha_k$ be the matrix defining  
$E^*$. Consider the fibre of $E^*$ at $a=[0,0,1]$ and use $x_1,x_2$ as affine coordinates. Then in the monad description  $E_a^*$ is defined as the quotient of $\C^3\otimes T^*(1)_a$ by the $3$-dimensional space spanned by   
$-A_{ij1}e_j\otimes dx_2+A_{ij2}e_j\otimes dx_1$ for $i=1,2,3$ where $e_i$ form a basis of $\C^3$.

The sections of $H(1)$ map $\PP(E)$ to the $5$-dimensional space 
$$\PP(H^0(\PP(E),H(1))^*)=\PP(H^0(\PP_2^*,E^*(1))^*)\subset \PP( \Hom(H^0(\PP^*_2,\Omega^1(2)),\C^3)).$$
The right hand side is an $8$-dimensional projective space of $3\times 3$ matrices $X$ and  the singular homomorphisms from  $H^0(\PP^*_2,\Omega^1(2))$ to $\C^3$ define a cubic determinantal  hypersurface $\det X =0$. 

\begin{prp} The linear system of the line bundle $H(1)$ maps $\PP(E)$ to the  intersection of the determinantal cubic hypersurface with 
$\PP(H^0(\PP(E),H(1))^*)$. 
\end{prp}
\begin{prf} 
 The section $\alpha_3=x_1dx_2-x_2dx_1$ of $T^*(2)$ vanishes at $a=[0,0,1]$ and so for any $v\in \C^3$, $v\otimes \alpha_3\in H^0(\PP_2^*,\C^3\otimes \Omega^1(2))$   maps to a section of $H^0(\PP_2^*,E^*(1))$ which vanishes at $a$. If $x\in \PP(E_a)$ then the evaluation map $\ev_x: H^0(\PP(E),H(1))\rightarrow H(1)_x$ defines the corresponding point for the projective embedding in $\PP(H^0(\PP(E),H(1))^*)$. 
 
 Under the map   $\PP(H^0(\PP(E),H(1))^*)\rightarrow \PP( \Hom(H^0(\PP^*_2,\Omega^1(2),\C^3))$ the point $x$ therefore  maps to a homomorphism for which $\alpha_3$ lies in the kernel, and hence a singular matrix.
\end{prf}

The cubic hypersurface $\det X=0$ has a singularity on the locus where the rank of $X$ is equal to one. The map $(v,w)\mapsto v\otimes w$ identifies this with $\PP_2\times\PP_2$. Since this is bilinear, the restriction of ${\mathcal O}(1)$ on $\PP_8$ to $\PP_2\times\PP_2$ is ${\mathcal O}(1,1)$, and so  its intersection with a generic $5$-dimensional space is a complete intersection of three sections. Since the canonical bundle of $\PP_2\times\PP_2$ is ${\mathcal O}(-3,-3)$, this means the intersection has trivial canonical bundle, i.e. is an elliptic curve of bidegree $(3,3)$. Recall now, that our Poisson deformations contain a distinguished copy of $D^{(2)}$, the zero set of the Poisson structure $\sigma_a$, and this is a $\PP_1$-bundle over the elliptic curve $\Pic^2(D)$. 

\begin{prp} The map $f:\PP(E)\rightarrow \PP_8$ defined by the line bundle  $H(1)$ collapses each fibre of $D^{(2)}\rightarrow \Pic^2(D)$ to a point and identifies $\Pic^2(D)$ with the singular locus of the cubic fourfold $f(\PP(E))$. \end{prp} 

\begin{prf} As remarked above, the fibres of  $\PP(T)\rightarrow \PP_2^*$, the exceptional divisor, collapse to points under the map and so $H(1)$ is trivial restricted to these. Choose a  point $z\in D\subset \PP_2$, then each line through $z$ meets $D$ again in a pair of points in the same divisor class, so we can identify $\PP(T_z)$ as a fibre  of  $D^{(2)}\rightarrow \Pic^2(D)$.   But this is a fibre of the exceptional divisor  -- the tangent directions at $z$. Hence $c_1(H(1))$ vanishes on this line.

But the deformations we constructed in Theorem \ref{def} preserved  the zero set of the Poisson structure and its complex structure so after deformation, $D^{(2)}$ is preserved and  the line bundle $H(1)$ for $\PP(E)$ is still trivial on a fibre of $D^{(2)}\rightarrow \Pic^2(D)$, which means that it is collapsed to a point under the linear system, and its image is the singular locus of the cubic fourfold. 
\end{prf} 

We have shown here how to recover the elliptic curve $D$, as an abstract curve, from the complex structure on the deformation of the Hilbert scheme -- it is the singular locus of the $-K/3$ model. Moreover it lies in $\PP_2\times\PP_2$ and so has two degree $3$ line bundles ${\mathcal O}(1,0)$ and ${\mathcal O}(0,1)$ on it. These provide the two parameters in the deformation -- the modulus of the curve and  a line bundle  ${\mathcal O}(1,-1)$. As we deform according to Theorem \ref{def} this line bundle changes linearly. 

In fact, if we embed an elliptic curve as a cubic in $\PP_2\times\PP_2$ in the standard way using theta functions we can describe explicitly the three linear equations in $3\times 3$ matrices defining $\PP_5\subset \PP_8$:
\begin{equation}
  \displaystyle\sum_{r\in\Z/3\Z}\frac{\theta_{j-i}(0)}
 {\theta_{r}(h)\theta_{j-i-r}(-h)}X_{{j-r},{i+r}}=0, \ \ i\neq j, i, j\in\Z/3\Z
 \end{equation}
and varying $h$ gives the deformation. We explain this in the next section.

\subsection{Sklyanin algebras}

The theta function formula above comes from the relations for an associative, non-commutative algebra due to Sklyanin \cite{Od} which is a deformation of a polynomial algebra in three variables. 
In \cite{Nev} Nevins and Stafford construct moduli spaces of modules over these Sklyanin algebras and show that the  moduli space of rank one torsion-free modules with $c_1=0$ and $\chi=1-n$ is a deformation of the Hilbert scheme $\PP_2^{[n]}$ with a natural  Poisson structure. They also give an explicit construction of this space as a quotient. 

More specifically, the algebra is defined by a 3-dimensional subspace of relations in $\C^3\otimes\C^3$ with the commutative polynomial algebra  defined by $\Lambda^2 \C^3\subset \C^3\otimes\C^3$. Invariantly, suppose $U,V,W$ are $3$-dimensional vector spaces with a homomorphism $U\rightarrow V\otimes W$, then this provides a $3\times 3$ matrix $Q\in V\otimes W$ of linear forms on $U$ whose determinant defines in general the equation of a cubic curve  $C\subset \PP(U)$. But on $C$, $Q$ is degenerate and the two maps $V^*\rightarrow W$ and $W^*\rightarrow V$ have kernels which define line bundles on $C$. Hence $C$ is embedded as a cubic curve in $\PP(U),\PP(V^*)$ and $\PP(W^*)$ (see e.g. \cite{Bea}).  This provides three line bundles $L_U,L_V,L_W$ with the relation (written additively) $L_V+L_W\cong 2L_U$. Up to equivalence $U,V,W$  define an elliptic curve and a translation $L_V-L_U$. This of course is the data that we have been experiencing in our deformation theory applied to $\PP_2^{[n]}$.

The construction in \cite{Nev} goes as follows: take elements $A\in \Hom(\C^n,\C^{2n+1})\otimes V$ and $B\in \Hom(\C^{2n+1},\C^n)\otimes W$ such that
$$BA\in \Hom(\C^n,\C^n)\otimes U\subset  \Hom(\C^n,\C^n)\otimes V\otimes W.$$

Then the moduli space of stable pairs $(A,B)$ with respect to the $GL(2n+1,\C)\times GL(n,\C)\times GL(n,\C)$ action is a deformation of the Hilbert scheme $\PP_2^{[n]}$.

  From Proposition \ref{h1} the generic deformation has a two-dimensional space of moduli, and indeed these two parameters  are the modulus of an elliptic curve together with a translation, so  this monad construction  applies to a generic deformation of the Hilbert scheme.

\vskip 1cm
 Mathematical Institute, 24-29 St Giles, Oxford OX1 3LB, UK
 
 hitchin@maths.ox.ac.uk

 \end{document}